\documentclass[12pt]{amsart}
\usepackage{amsmath,amsfonts,amsthm,amsopn,amssymb,mathtools,stmaryrd}
\usepackage{cite,marginnote}
\usepackage{amsmath}
\usepackage{color,enumitem,graphicx}
\usepackage[colorlinks=true,urlcolor=blue,
citecolor=red,linkcolor=blue,linktocpage,pdfpagelabels,
bookmarksnumbered,bookmarksopen]{hyperref}
\usepackage[english]{babel}

\usepackage[left=2.75cm,right=2.75cm,top=3cm,bottom=3cm]{geometry}

\usepackage[hyperpageref]{backref}

\numberwithin{equation}{section}

\makeindex
\newtheorem{theorem}{Theorem}[section]
\newtheorem{proposition}[theorem]{Proposition}
\newtheorem{lemma}[theorem]{Lemma}
\newtheorem{remark}[theorem]{Remark}
\newtheorem{example}[theorem]{Example}
\newtheorem{corollary}[theorem]{Corollary}
\newtheorem{definition}[theorem]{Definition}
\newtheorem{theoremletter}{Theorem}
\newtheorem{lemmaletter}{Lemma}

\newcommand{\ud}{\mathrm{d}}
\newcommand{\RN}{\mathbb R^N}
\newcommand{\om}{\Omega}

\newcommand{\iy}{\infty}

\newcommand{\s}{\section}
\newcommand{\dd}{\delta}
\newcommand{\DD}{\Delta}
\newcommand{\g}{\gamma}

\newcommand{\na}{\nabla}

\newcommand{\la}{\lambda}

\newcommand{\pa}{\partial}

\newcommand{\R}{\mathbb R}
\newcommand{\al}{\alpha}

\newcommand{\ti}{\tilde}

\newcommand{\re}[1]{(\ref{#1})}

\newcommand{\rg}{\rightarrow}

\newcommand{\lan}{\langle}
\newcommand{\ran}{\rangle}
\newcommand{\e}{\varepsilon}
\newcommand{\vp}{\varphi}

\newcommand{\lab}{\label}
\newcommand{\bt}{\begin{theorem}}
\newcommand{\et}{\end{theorem}}
\newcommand{\bl}{\begin{lemma}}
\newcommand{\el}{\end{lemma}}
\newcommand{\bd}{\begin{definition}}
\newcommand{\ed}{\end{definition}}
\newcommand{\bc}{\begin{corollary}}
\newcommand{\ec}{\end{corollary}}
\newcommand{\bp}{\begin{proof}}
\newcommand{\ep}{\end{proof}}
\newcommand{\bx}{\begin{example}}
\newcommand{\ex}{\end{example}}
\newcommand{\bi}{\begin{exercise}}
\newcommand{\ei}{\end{exercise}}
\newcommand{\bo}{\begin{proposition}}
\newcommand{\eo}{\end{proposition}}
\newcommand{\br}{\begin{remark}}
\newcommand{\er}{\end{remark}}
\newcommand{\be}{\begin{equation}}
\newcommand{\ee}{\end{equation}}
\newcommand{\ba}{\begin{align}}
\newcommand{\ea}{\end{align}}
\newcommand{\bn}{\begin{enumerate}}
\newcommand{\en}{\end{enumerate}}
\newcommand{\bg}{\begin{align*}}
\newcommand{\bcs}{\begin{cases}}
\newcommand{\ecs}{\end{cases}}

\newcommand{\pl}{\partial}
\newcommand{\vr}{\varepsilon}

\newcommand{\bean}{\begin{eqnarray*}}
\newcommand{\eean}{\end{eqnarray*}}
\newcommand{\loc}{\operatorname{\rm loc}}

\title[A priori estimates and positivity for semiclassical ground states]{A priori estimates and positivity for semiclassical ground states for systems of critical Schr\"{o}dinger equations in dimension two}

\author[D.~Cassani]{Daniele Cassani$^\text{1}$}
\author[J.~Zhang]{Jianjun Zhang$^\text{2}$}

\address[D. Cassani]{\newline\indent Dip. di Scienza e Alta Tecnologia
\newline\indent
Universit\`{a} degli Studi dell'Insubria\\
\newline\indent via Valleggio 11, 22100 Como,Italy
\newline\indent and
\newline\indent RISM--Riemann International School of Mathematics\\
\newline\indent Villa Toeplitz, Via G.B. Vico, 46 -- 21100 Varese}
\email{\href{mailto:Daniele.Cassani@uninsubria.it}{Daniele.Cassani@uninsubria.it}}

\address[J.~Zhang]{\newline\indent College of Mathematics and Statistics
\newline\indent
Chongqing Jiaotong University
\newline\indent
Chongqing 400074, PR China
\newline\indent and
\newline\indent Dip. di Scienza e Alta Tecnologia
\newline\indent
Universit\`{a} degli Studi dell'Insubria
\newline\indent
via Valleggio 11, 22100 Como,Italy}
\email{\href{mailto:zhangjianjun09@tsinghua.org.cn}{zhangjianjun09@tsinghua.org.cn}}

\thanks{(1) Corresponding author: \texttt{daniele.cassani@uninsubria.it}}
\thanks{(2) The second named author was partially supported by the Science Foundation of Chongqing Jiaotong University(15JDKJC-B033)}


\begin{document}

\begin{abstract}
We consider in the whole plane the following Hamiltonian coupling of Schr\"{o}dinger equations
\begin{equation*}\left\{
\begin{array}{ll}
-\DD u+V_0u=g(v)\\
-\DD v+V_0v=f(u)\\
\end{array}
\right. \end{equation*}
where $V_0>0$, $f,g$ have critical growth in the sense of Moser. We prove that the (nonempty) set $S$ of ground state solutions is compact in $H^1(\R^2)\times H^1(\R^2)$ up to translations. Moreover, for each $(u,v)\in S$, one has that $u,v$ are uniformly bounded in $L^\iy(\R^2)$ and uniformly decaying at infinity. Then we prove that actually the ground state is positive and radially symmetric. We apply those results to prove the existence of semiclassical ground states solutions to the singularly perturbed system
\begin{equation*}
\begin{cases}
-\vr^2\Delta \varphi+V(x)\varphi=g(\psi)&\\
-\vr^2\Delta \psi+V(x)\psi=f(\varphi)
\end{cases}
\end{equation*}
where $V\in \mathcal{C}(\mathbb{R}^2)$ is a Schr\"{o}dinger potential bounded away from zero. Namely, as the adimensionalized Planck constant $\vr\to 0$, we prove the existence of minimal energy solutions which concentrate around the closest local minima of the potential with some precise asymptotic rate. 
\end{abstract}
\maketitle

\s{Introduction}
\noindent
Consider in the whole $\mathbb{R}^2$ the following system of coupled Schr\"odinger equations
\be\lab{q1} \left\{
\begin{array}{ll}
\displaystyle -\e^2\DD \vp+V(x)\vp=\frac{\partial H(\varphi,\psi)}{\partial\psi}\\
\displaystyle -\e^2\DD \psi+V(x)\psi=-\frac{\partial H(\varphi,\psi)}{\partial\varphi}
\end{array}
\right. 
\ee
where $\e>0$, the external Schr\"odinger potential $V\in C(\R^2,\R)$ enjoys the following condition:  
\begin{itemize}
\item [$(V)$] $0<V_0:=\inf_{\R^2}V(x)<\lim_{|x|\rg\iy}V(x)=V_\iy\leq\infty$.
\end{itemize}
The Hamiltonian has the following form $H(\varphi,\psi)=G(\psi)-F(\varphi)$, with $F(t)=\int_0^t f(\tau)\, \ud \tau$ and $G(t)=\int_0^t g(\tau)\, \ud \tau$ and the nonlinearities $f,g\in C(\R,\R)$  satisfy the following hypotheses:
\begin{itemize}
\item [$(H1)$] $f(t)=o(t)$ and $g(t)=o(t)$, as $t\rg0$;
\item [$(H2)$] There exists $\theta>2$ such that for any $t\not=0$,
\begin{center}
$0<\theta F(t)\le f(t)t$ and $0<\theta G(t)\le g(t)t$;
\end{center}
\item [$(H3)$] There exists $M>0$ such that for any $t\not=0$,
\begin{center}
$0<F(t)\le Mf(t)$ and $0<G(t)\le Mg(t)$;
\end{center}
\item [$(H4)$] $f(t)/|t|$ and $g(t)/|t|$ are strictly increasing for $t\not=0$.
\end{itemize}

\noindent As a consequence of the Pohozaev-Trudinger-Moser inequality for which the Sobolev space $H^1$ embeds into the space of functions such that $e^{\alpha u^2}\in L^1$, the following notion of critical growth in dimension two was first introduced in \cite{AY,DMR} (in the case of bounded domains):
\bd\label{cgdef} A function $f:\R^+\rightarrow\R^+$ has critical growth in the sense of Pohozaev-Trudinger-Moser inequality, if there exists $\al_0>0$ such that
\begin{align*}
\lim_{t\to + \infty} \frac {f(t)}{e^{\alpha t^2}}=
\begin{cases}
0 & \text{if } \al > \al_0
\\
+ \infty & \text{if } \al < \al_0
\end{cases}
\end{align*}
\ed

\noindent It will be crucial in what follows the following growth assumptions:
\medskip
\begin{itemize}
\item [$(H5)$] $\displaystyle\liminf_{|t|\rg\iy}\frac{tf(t)}{e^{\al_0t^2}}\ge\beta_0>\frac{2e}{\al_0}V_0$ and $\displaystyle\liminf_{|t|\rg\iy}\frac{tg(t)}{e^{\al_0t^2}}\ge\beta_0>\frac{2e}{\al_0}V_0.$
\end{itemize}
\medskip
\noindent It is well known, both from the theoretical point of view as well as from that of applications, that minimal energy solutions, the so-called ground states, play a fundamental role, see e.g. \cite{BF}. In what follows we will focus on this class of solutions. In particular, to investigate the sign of ground state solutions to \re{q1}, we require in addition the following condition:
\begin{itemize}


\item [$(H6)$] There exist $p, q>1$ such that $f(t)\geq t^q$ and $g(t)\geq t^p$ for small $t>0$;
\end{itemize}
\noindent As a reference model take $F(t)=|t|^p (\mathrm{e}^{4\pi t^2}-1)$ and $G(t)=|t|^q(\mathrm{e}^{4\pi t^2}-1)$ with $p,q>2$ and $\al_0=4\pi$ which clearly satisfy $(H1)$-$(H6)$.

Our main result reads as follows:
\bt\lab{Th1} Assume condition $(V)$ and that $f,g$ have critical growth in the sense of Definition \ref{cgdef} and satisfy $(H1)$--$(H5)$. Then, for sufficiently small $\e>0$, $(\ref{q1})$ admits a least energy solution $z_\e=(\vp_\e,\psi_\e)\in H^1(\R^2)\times H^1(\R^2)$. Moreover, the following properties hold:
\begin{itemize}
\item [$(i)$] let $x_\e^1, x_\e^2, x_\e$ be any maximum point of $|\vp_\e|, |\psi_\e|, |\vp_\e|+|\psi_\e|$ respectively, then, setting $$\mathcal{M}\equiv\{x\in \R^2: V(x)=V_0\}$$ one has
$$
\lim_{\e\rg
    0}\mbox{dist}(x_\e,\mathcal{M})=0\:\text{ and }\: \lim_{\e\rg0}|x_\e^i-x_\e|=0,\quad i=1,2.
$$
Furthermore, $(\vp_\e(\e x+x_\e),\psi_\e(\e x+x_\e))$ converges (up to a subsequence) as $\e\to 0$ to a ground state solution of
\begin{align*}
\left\{
\begin{array}{ll}
-\DD u+V_0u=g(v)\\
-\DD v+V_0v=f(u)
\end{array}
\right.
\end{align*}
\item [$(ii)$] if in addition $(H6)$ holds, then replacing $f$ and $g$ above with their odd extensions, for $\e>0$ small enough, up to changing sign $u_\e, v_\e>0$ in $\R^2$ and $x_\e^1,x_\e^2$ are the unique global maximum points of $u_\e,v_\e$ respectively and which also enjoy the following
$$\lim_{\e\rg0}|x_\e^1-x_\e^2|/\e=0.$$
Moreover, for some $c,C>0$ one has $$|\vp_\e(x)|\le C\exp(-\frac{c}{\e}|x-x_\e^1|),\,\,|\psi_\e(x)|\le C\exp(-\frac{c}{\e}|x-x_\e^2|), \,\, x\in\R^2;$$
\end{itemize}
\et
(Without loss of generality, throughout the paper we may assume $0\in\mathcal{M}$.)
\br\label{remark1} Let us point out a few comments on the conditions we assume in Theorem \ref{Th1}: 
\begin{itemize}
\item Actually the Ambrosetti-Rabinowitz condition $(H2)$ can be replaced by the following slightly weaker assumption:
\begin{itemize}
\item [$(H2)'$] There exists $\theta>2$ such that for any $t\not=0$,
\begin{center}
$0<2F(t)\le f(t)t$ and $0<\theta G(t)\le g(t)t$,
\end{center}
or equivalently
\begin{center}
$0<\theta F(t)\le f(t)t$ and $0<2G(t)\le g(t)t$.
\end{center}
\end{itemize}
\item We also point out that conditions $(H2)$ and $(H4)$ are weaker than the following assumption:
\begin{itemize}
\item [$(H)$] $f,g\in C^1(\R,\R)$ and there exists $\dd'>0$ such that for any $s\not=0$,
$$
0<(1+\dd')f(s)s\le f'(s)s^2\,\,\mbox{and}\,\,\, 0<(1+\dd')g(s)s\le g'(s)s^2
$$
which appears in the literature, see \cite{dsr,Pisto,Ramos1}.
\end{itemize}
\item Hypotheses $(H6)$ and $(H)$ can be also found in \cite{dsr}. Clearly hypothesis $(H)$ yields $sf(s)\le f(1)|s|^{2+\dd'}$ and $sg(s)\le g(1)|s|^{2+\dd'}$ if $|s|\le1$. Let us point out that in the present paper we do not require $sf(s),sg(s)$ to be less than $|s|^r$ near the origin for some $r>2$.
\end{itemize}
\er
Systems of the form \eqref{q1} have been largely investigated in the last three decades being a prototype in many different applications, where they model for instance the minimal energy interaction between nonlinear fields, see \cite{BF,yang}. The scenario changes remarkably from the higher dimensional case $N\geq 3$ to the planar case $N=2$. In particular, $N=2$ affects the notion of critical growth which is the maximal admissible growth for the nonlinearities in order to preserve the variational structure of the problem; we refer to \cite{dcassani2,CST,CT} for a discussion on this topic and to \cite{bst,ruf} for a survey on systems of the form \eqref{q1} in the case of bounded domains. As far as we are concerned with minimal energy solutions in the whole space, existence results have been first established in \cite{boyan}, see also \cite{Weth}, in the higher dimensional case and then recently extended to $N=2$ in \cite{DJJ} , where the Trudinger-Moser critical case is covered, see also \cite{bsrt,dsr}. Qualitative properties of minimal energy solutions such as symmetry and positivity have been investigated in the higher dimensional case in \cite{Sirakov1,syso,dsr}, see also \cite{dgp,QS} for closely related results. Always in dimension $N\geq 3$, a priori estimates have been obtained in \cite{DY}. A priori bounds open the way to investigate the existence and concentrating behavior, as $\e\to 0$, of the so-called semiclassical states. From the point of view of Physics, these solutions live on the interface between quantum and classical Mechanics, in the sense that the field behaves like a Newtonian particle as $\e\to 0$, see \cite{EvansZ} for a survey on the topic and references therein. Semiclassical states for singularly perturbed Schr\"{o}dinger systems have been studied in the higher dimensional case in \cite{ASY,Ramos1,dlz}. 
   
\noindent Finally, let us mention that the question weather the ground state we find is unique, seems to be out of reach at the moment. This is still a challenging open problem even in the subcritical case as well as in higher dimensions in which uniqueness of positive solutions (not necessarily ground states) is known just in a few particular cases such as Lane-Emden systems \cite{rd}. More in general, the matter of uniqueness of ground states, even in cases in which one has multiplicity of positive solutions, remains open even for the single equation.  

\subsection*{Overview} The paper is organized as follows: in Section \ref{lp} we begin with studying a limit problem for system \eqref{q1}. Here we complete the work initiated in \cite{DJJ}, where the existence of ground states is proved, by establishing a priori estimates, regularity, symmetry and qualitative properties of solutions. Here we use a suitable Nehari manifold approach in the spirit of Pankov \cite{Pankov} combined with Moser type techniques, as everything is set in dimension two and in presence of Moser critical growth. In particular we exploit those preliminary results to prove positivity of ground states solutions in a quite general setting, as developed throughout Section \ref{sign_s}. Then, Section \ref{semiclassical_s}  is devoted to apply the informations previously obtained on the limit problem, to analyze the concentrating behavior of semiclassical solutions from the point of view of localizing bumps as well as of deriving the asymptotic rate of concentration. Here the presence of critical Moser's growth requires some delicate energy estimates which we then apply to establish compactness. 
\s{The limit problem}\label{lp}

\noindent
By denoting $u_\e(x)=\vp(\e x),v_\e(x)=\psi(\e x)$ and $V_\e(x)=V(\e x)$, \re{q1} is equivalent to
\begin{align*} \left\{
\begin{array}{ll}
-\DD u_\e+V_\e(x)u_\e=g(v_\e)\\
-\DD v_\e+V_\e(x)v_\e=f(u_\e)
\end{array}
\right. \end{align*}
in the whole plane. Let $x_0\in\R^2$ and assume $u_\e(\cdot+\frac{x_0}{\e})\rg u$, $v_\e(\cdot+\frac{x_0}{\e})\rg v$ in
$C^1_{loc}(\R^2)$, if $V_0=V(x_0)$ then one has 
\be\lab{q11} \left\{
\begin{array}{ll}
-\DD u+V_0u=g(v)\\
-\DD v+V_0v=f(u)
\end{array}
\right. \ee
which is the so-called limit problem of \re{q1}. Recently, D.\ G. De\ Figueiredo, J. M. do \'O and J. Zhang established in \cite{DJJ} the existence of ground state solutions to \re{q11}, precisely
\begin{theoremletter}\lab{a}
{\rm (Theorem 1.3 in \cite{DJJ})} {\it Suppose that $f,g$ have critical growth and satisfy $(H1)$--$(H5)$. Then \eqref{q11} admits a ground state solution $(u,v)\in H^1(\R^2)\times H^1(\R^2)$.}
\end{theoremletter}
\noindent Denote by $\mathcal{S}$ the set of of ground state solutions to system \re{q11}, then by Theorem A $\mathcal{S}\not=\emptyset$. Here we investigate the regularity and qualitative properties of the ground state solutions to \re{q11}. Precisely, we prove the following results:

\bt\lab{Th2} Suppose $f,g$ have critical growth and satisfy $(H1)$--$(H5)$. Then the following hold true:
\begin{itemize}
\item [$(i)$] $(u,v)\in \mathcal{S}\Rightarrow u,v\in L^{\iy}(\R^2)\cap C_{loc}^{1,\g}(\R^2)$ for some $\g\in(0,1)$;
\item [$(ii)$] let $x_z\in\R^2$ be the maximum point of $|u(x)|+|v(x)|$, then the set $$\{(u(\cdot+x_z),v(\cdot+x_z))\,|\, (u,v)\in \mathcal{S}\}$$ is compact in $H^1(\R^2)\times H^1(\R^2)$;
\item [$(iii)$] $0<\inf\{\|u\|_{\infty},\|v\|_{\infty}\,|\, (u,v)\in \mathcal{S}\}\le \sup\{\|u\|_{\infty},\|v\|_{\infty}\,|\, (u,v)\in \mathcal{S}\}<\iy$;
\item [$(iv)$] $u(x+x_z)\rg 0$ and $v(x+x_z)\rg 0$, as $|x|\rg\iy$ uniformly for any $z=(u,v)\in \mathcal{S}$, where $x_z$ is given in $(ii)$;
\item [$(v)$] for any $(u,v)\in \mathcal{S}$, the following Pohozaev-type identity holds
$$
\int_{\R^2}(F(u)+G(v)-V_0uv)\,\ud x=0.
$$
\end{itemize}
\et
\bt\lab{sign} Assume in addition to the hypotheses of Theorem \ref{Th2} that also $(H6)$ holds. Then, replacing $f$ and $g$ in Theorem \ref{Th2} with their odd extensions, for any $(u,v)\in\mathcal{S}$ one has $u,v\in C^2(\R^2)$ and $uv>0$ in $\R^2$. Moreover, there exists some point
$x_0\in\R^2$ such that $u,v$ are radially symmetric with respect to the same point $x_0$, namely $u(x)=u(|x-x_0|)$, $v(x)=v(|x-x_0|)$ and setting $r=|x-x_0|$, one has for $r>0$
$$
\frac{\pa u}{\pa r}<0\quad \text{ and }\quad \frac{\pa v}{\pa r}<0
$$
as well as 
$$\DD u(x_0)<0\quad  \text{ and }\quad \DD v(x_0)<0.$$
Moreover, there exist $C,c>0$, independent of $z=(u,v)\in \mathcal{S}$, such that $$|D^{\al}u(x)|+|D^{\al}v(x)|\le C\exp(-c|x-x_0|),\quad x\in \R^2,\,|\al|=0,1$$
\et


\subsection{The functional setting: a generalized Nehari manifold}

\renewcommand{\theequation}{2.\arabic{equation}}

Let $H^1(\R^2)$ be the usual Sobolev space endowed with the inner product
$$
(u,v)_{H^1}:=\int_{\R^2}\na u\na v+V_0uv,\ \ \|u\|_{H^1}^2:=(u,u)_{H^1}, u,v\in H^1(\R^2).
$$
and set $E=H^1(\R^2)\times H^1(\R^2)$ with the inner product
$$
(z_1,z_2):=(u_1,u_2)_{H^1}+(v_1,v_2)_{H^1},\ \ z_i=(u_i,v_i)\in E, i=1,2.
$$
Clearly we have the space decomposition $E=E^+\oplus E^-$, where
$$
E^+:=\{(u,u)\,|\, u\in H^1(\R^2)\}\ \ \ \mbox{and}\ \ \ E^-:=\{(u,-u)\,|\, u\in  H^1(\R^2)\}.
$$
For each $z=(u,v)\in E$, one has $$z=z^++z^-=((u+v)/2,(u+v)/2)+((u-v)/2,(v-u)/2).$$

\vskip0.1in
Weak solutions to \eqref{q11} are the critical points of the associated energy functional
$$
\Phi(z):=\int_{\R^2}\na u\na v+V_0uv-I(z),\ \ z=(u,v)\in E,
$$
where $I(z)=\int_{\R^2}F(u)+G(v)$.
Using the above notation we have
\be\lab{y1}
\Phi(z):=\frac{1}{2}\|z^+\|^2-\frac{1}{2}\|z^-\|^2-I(z),
\ee
which emphasizes the strongly indefinite nature of $\Phi$ which however, by the hypotheses on $f$ and $g$, is of class $C^1(E,\R)$ and
\be\lab{y2}
I(0)=0, \ \lan I'(z), z\ran>2I(z)>0,\ \ \mbox{for all}\ \ z\in E\setminus\{0\}.
\ee
On one hand, if $z=(u,v)\in E\setminus \{0\}$ such that $\Phi'(z)=0$, then by $(H2)$
\begin{equation*}
\Phi(z)=\Phi(z)-\frac{1}{2}\lan \Phi'(z),z\ran=\int_{\R^2}\frac{1}{2}f(u)u-F(u)+\frac{1}{2}g(u)u-G(u)>0.
\end{equation*}
On the other hand, if $z=(u,-u)\in E^-$, we have by $(H2)$ 
$$
\Phi(z)=-\int_{\R^2}(|\na u|^2+V_0u^2)-\int_{\R^2} F(u)+G(-u)\le 0.
$$
As a consequence, if $z\in E$ is a nontrivial critical point of $\Phi$, then necessarily $z\in E\setminus E^-$. This motivates the introduction of the following generalized Nehari manifold, due to Pankov \cite{Pankov} and then used also in \cite{Szulkin, Weth, DJJ}:
$$
\mathcal{N}:=\{z\in E\setminus E^-: \lan \Phi'(z),z\ran=0, \lan \Phi'(z),\vp\ran=0\ \mbox{for all}\ \ \vp\in E^-\}.
$$
Let
$$
c_\ast:=\inf_{z\in\mathcal{N}}\Phi(z)
$$
then $c_\ast$ is called the least energy level of system \re{q11}. In \cite{DJJ} the authors proved that $c_\ast\in(0,4\pi/{\al_0})$ and that it is achieved on $\mathcal{N}$.


\subsection{Proof of Theorem \ref{Th2}}
\renewcommand{\theequation}{3.\arabic{equation}}
Let $\{z_n\}\subset S$, namely 
\be\lab{pss}
\Phi(z_n)=c_\ast \ \ \mbox{and}\ \ \Phi'(z_n)=0, \quad \forall n\in\mathbb{N}
\ee
We carry out the proof of $(ii)$ of Theorem \ref{Th2} through the following four steps:
\begin{itemize}
\item We first prove that $\{z_n\}$ is bounded in $E$ (Proposition \ref{o1});
\item In Proposition \ref{nv} we prove that there exisst $\{y_n\}\subset\R^2$ and $z_0\not={\bf 0}$ such that $z_n(\cdot+y_n)\rightharpoonup z_0$ in $E$ and $z_n(\cdot+y_n)\xrightarrow{a.e.}z_0$ in
$\R^2$, as $n\rg\iy$;
\item In Proposition \ref{o11} we show that $z_0$ is actually a critical point of $\Phi$;
\item Finally in Proposition \ref{con} we prove that $z_0\in \mathcal{S}$ and that actually $z_n(\cdot+y_n)\longrightarrow z_0$ strongly in $E$, as $n\to \infty$.
\end{itemize}
In the proof of the Proposition \ref{o1} below we will use the following lemma which we borrow from \cite{Fi}:
\begin{lemmaletter}\label{RUF1} {\it The following inequality holds
\[
    s\text{ }t\leq \left\{
\begin{array}{ll}
(e^{t^{2}}-1)+s(\log s)^{1/2}, \; & \text{ for
all }t\geq 0 \text{ and }s\geq e^{1/4}; \\
(e^{t^{2}}-1)+\frac{1}{2}s^{2}, \; & \text{ for all } t \geq 0
\text{ and }0 \leq s\leq e^{1/4}.
\end{array}
\right.
\]}
\end{lemmaletter}
\noindent The proofs of Proposition \ref{o1} and \ref{o11} are similar to \cite{DJJ}, however for the sake of completeness we give the details.
\bo\lab{o1} There exists $C>0$ such that for all $n\in\mathbb{N}$:
\begin{itemize}

\item [$1)$] $\|z_n\|=\|(u_n,v_n)\|\le C$;

\item [$2)$] $\int_{\R^2}f(u_n)u_n\, \ud x\le C$ and $\int_{\R^2}g(v_n)v_n\, \ud x\le C$;

\item [$3)$] $\int_{\R^2}F(u_n)\, \ud x\le C$ and $\int_{\R^2}G(v_n)\, \ud x\le C$.
\end{itemize}
\eo
\bp
From $\lan\Phi'(z_n),z_n\ran=0$ we have
\begin{equation}\label{pitomba}
 2 \int_{\R^2} (\nabla u_n \nabla v_n+V_0u_nv_n)\, \ud x - \int_{\R^2}
f(u_n)u_n\, \ud x -  \int_{\R^2} g(v_n)v_n \, \ud x =0.
\end{equation}
Recalling that
$$
\Phi(z_n)=\int_{\R^2}(\nabla u_n \nabla v_n+V_0u_nv_n)\, \ud x-\int_{\R^2}(F(u_n)+G(v_n))\, \ud x=c_\ast
$$
we obtain by $(H_3)$ the following 
\begin{align*}
\int_{\R^2} [f(u_n)u_n + g(v_n)v_n]\, \ud x &= 2\int_{\R^2} [F(u_n) +  G(v_n)]\, \ud x + 2c_\ast \\
& \leq \frac{2}{\theta}\int_{\R^2} [f(u_n)u_n + g(v_n)v_n]\, \ud x + 2c_\ast.
\end{align*}
Thus
\begin{equation}\label{goiaba}
\int_{\R^2} [f(u_n)u_n + g(v_n)v_n]\, \ud x  \leq  \frac{2c_\ast\theta}{\theta-2}.
\end{equation}
From $ \lan\Phi'(z_n),(v_n,0)\ran=0$ and $ \lan\Phi'(z_n),(0,u_n)\ran=0$, we have
$$
\| v_n \|^2-\int_{\R^2}
f(u_n)v_n\, \ud x=0\ \ \mbox{and}\ \ \| u_n \|^2-\int_{\R^2}
g(v_n)u_n\, \ud x=0.
$$
Let $U_n=u_n/ \| u_n\|$ and  $ V_n = v_n / \| v_n \| $, then
\begin{align}
\| v_n \|   &  =  \int_{\R^2}
f(u_n)V_n\, \ud x   \label{laranja},\\
\| u_n \|   &  =  \int_{\R^2}
g(v_n)U_n\, \ud x  . \label{limao}
\end{align}
By $(H1)$, there exist $\beta>0$ and $C_\beta>0$ such that
$$
f(t)\le C_\beta e^{\beta t^2}\quad \text{ and }\quad g(t)\le C_\beta e^{\beta t^2}\ \ \mbox{for all} \ \ t\ge0.
$$
Moreover, there exists $C_1>0$ such that for all $n$
$$
f(u_n(x))\le C_1 u_n(x)\ \ \mbox{for}\ \ x\in\{\R^2 : f(u_n(x))/C_\beta \leq e^{1/4}\}.
$$
Setting $ t = V_n $ and $ s =
f(u_n)/C_\beta $ in Lemma \ref{RUF1} then by $(H1)$-$(H2)$ together with the Pohozaev-Trudinger-Moser
inequality, we get 
\begin{align*}
\int_{\R^2} f(u_n)V_n\, \ud x& \leq  C_\beta \int_{\{x \in \R^2 : f(u_n(x))/C_\beta \geq e^{1/4} \}}
\frac{1}{C_\beta}f(u_n) \left[\log (\frac{1}{C_\beta} f(u_n))\right]^{1/2}\, \ud x \\
&+ \frac{1}{2}\int_{\{x \in \R^2 : f(u_n(x))/C_\beta \leq e^{1/4} \}}
\frac{1}{C_\beta} (f(u_n))^2\, \ud x+C_\beta \int_{\R^2} (e^{V_n^{2}}-1)
\, \ud x \\
& \leq C_2 + (\beta^{1/2}+C_1/(2C_\beta) \int_{ \R^2 }  f(u_n) u_n \, \ud x,
\end{align*}
for some constant $ C_2>0$. This estimate together with
(\ref{laranja}) implies, for some constant $ c_1 > 0 $, that
\begin{equation}\label{maravilha}
\| v_n \|  \leq  c_1(1 + \int_{ \R^2 }  f(u_n) u_n\, \ud x)
\end{equation}
and similarly
\begin{equation}\label{macacheira}
\| u_n \|  \leq  c_1(1 + \int_{ \R^2 }  g(v_n) v_n\, \ud x).
\end{equation}
From \re{maravilha}, \re{macacheira} and \re{goiaba} it follows the first claim $1)$. Then, by (\ref{goiaba}) and $(H_3) $ we
obtain the remaining bounds $2)$ and $3)$. 
\ep
Next we prove that, up to translations, $\{z_n\}$ has a nontrivial weak limit. Clearly $(u_n,v_n)$ satisfies just one of the following conditions:
\begin{itemize}

\item[] ({\it Vanishing})
$
\quad \lim_{n\rg\iy}\sup_{y\in\R^2}\int_{B_R(y)}(u_n^2+v_n^2)\, \ud x=0\ \ \mbox{for all}\ \ R>0;
$

\vspace{0,1cm}
\item[] ({\it Nonvanishing}) there exist $\nu>0$, $R_0>0$ and $\{y_n\}\subset\R^2$ such that
$$
\lim_{n\rg\iy}\int_{B_{R_0}(y_n)}(u_n^2+v_n^2)\, \ud x\ge\nu.
$$
\end{itemize}
We borrow from \cite{Fi} the following lemma:
\begin{lemmaletter}\lab{l3.3}{\it
Let $\Omega\subset\R^2$ be a bounded domain and $f\in C(\R,\R)$. Let $\{u_n\}\subset L^1(\Omega)$ be such that $u_n\rg u$ strongly in $L^1(\Omega)$,
$$
f(u_n)\in L^1(\Omega)\ \ \mbox{and}\ \ \ \int_{\Omega}|f(u_n)u_n|\, \ud x\le C, n\ge1
$$
for some $C>0$. Then, up to a subsequence we have 
$$
f(u_n)\rg f(u)\ \ \mbox{strongly in}\ \ L^1(\Omega)\ \ \mbox{as}\ \ n\rg\iy.
$$}
\end{lemmaletter}

\bo\lab{nv}
Vanishing does not occur.
\eo
\bp We know from \cite{DJJ} that $c_\ast\in(0,4\pi/\al_0)$, hence for some $\dd>0$ sufficiently small one has $c_\ast\in(0,4\pi/\al_0-\dd)$.
Assume by contradiction that vanishing occurs, namely 
$$
\lim_{n\rg\iy}\sup_{y\in\R^2}\int_{B_R(y)}(u_n^2+v_n^2)\, \ud x=0\ \ \mbox{for all}\ \ R>0,
$$
then Lions's lemma \cite{lionslemma} yields $u_n\rg0, v_n\rg0$ strongly in $L^s(\R^2)$ for any $s>2$.

Let us divide the proof into two steps:
\vskip0.1in
{\bf Step 1.} We claim 
$$
\lim_{n\rg\iy}\int_{\R^2}F(u_n)\, \ud x=0\ \ \mbox{and}\ \ \lim_{n\rg\iy}\int_{\R^2}G(v_n)\, \ud x=0.
$$
Indeed, by Lemma \ref{l3.3}, for any $R>0$ one has $f(u_n)\rg0$ and $g(v_n)\rg0$ strongly in $L^1(B_R(0))$ as $n\rg\iy$. Then by $(H3)$ and the Lebesgue dominated convergence theorem,
\be\lab{y4.1}
\lim_{n\rg\iy}\int_{B_R(0)}F(u_n)\, \ud x=0\ \ \mbox{and}\ \ \lim_{n\rg\iy}\int_{B_R(0)}G(v_n)\, \ud x=0.
\ee
In order to prove the claim, it is enough to prove that for any $\e>0$, there exists $R>0$ such that for $n$ large enough,
\be\lab{y4.2}
\int_{\R^2\setminus B_R(0)}F(u_n)\, \ud x\le\e\quad\text{and}\quad \int_{\R^2\setminus B_R(0)}G(v_n)\, \ud x\le\e.
\ee
By $(H3)$ and Proposition \ref{o1}, for any $K>0$ and $n$,
$$
\int_{\{x\in\R^2\setminus B_R(0): |u_n(x)|\ge K\}}F(u_n)\le\frac{M}{K}\int_{\{x\in\R^2\setminus B_R(0): |u_n(x)|\ge K\}}f(u_n)u_n\le\frac{MC}{K}.
$$
Then choosing $K>0$ large enough, we get that for all $n$
\be\lab{y4.3}
\int_{\{x\in\R^2\setminus B_R(0): |u_n(x)|\ge K\}}F(u_n)\le\frac{\e}{2}.
\ee
By $(H1)$, for any $\rho>0$ there exists $C_{\rho,K}>0$ such that
$$
F(t)\le \rho t^2+C_{\rho,K}t^4,\ \ \ |t|\le K.
$$
Recalling that $u_n\rg 0$ strongly in $L^4(\R^2)$, we obtain
$$
\limsup_{n\rg\iy}\int_{\{x\in\R^2\setminus B_R(0): |u_n(x)|\le K\}}F(u_n)\le\rho\sup_{n}\|u_n\|_2^2.
$$
By Proposition \ref{o1} and since $\rho$ is arbitrary, for $n$ large enough we get 
\be\lab{y4.4}
\int_{\{x\in\R^2\setminus B_R(0): |u_n(x)|\le K\}}F(u_n)\le\frac{\e}{2}.
\ee
Thus \re{y4.3} and \re{y4.4} yield the first bound in \re{y4.2} and similarly one gets the second bound.

\vskip0.1in
{\bf Step 2.} We claim that $c_\ast=0$, from which the contradiction follows as we know $c_*>0$. We need the following inequality used in \cite[Lemma 4.1]{Souza}
\be\lab{inequa}
t\ s\le t^2(e^{t^2}-1)+s(\log{s})^{\frac{1}{2}},\ \ \mbox{for all}\ \ (t,s)\in[0,\iy)\times[e^{\frac{1}{\sqrt[3]{4}}},\iy).
\ee
By Step 1,
\be\lab{y4.5}
\lim_{n\rg\iy}\int_{\R^2}(\na u_n\na v_n+V_0u_nv_n)=c_*\,.
\ee
If $u_n\rg 0$ or $v_n\rg0$ strongly in $H^1(\R^2)$ as $n\rg\iy$, then \re{y4.5} directly yields $c_*=0$. Therefore, let us assume $\inf_{n\ge1}\|u_n\| \ge b > 0$. Note that
\begin{equation}\label{unvn4}
\| u_n \|^2 = \int_{\R^2} g(v_n) u_n \, \ud x.
\end{equation}
By $(H1)$, for any fixed $\e>0$, there exists $C_\e>0$ such that
$$
f(t), g(t)\le C_\e e^{(\al_0+\e)t^2}\ \ \mbox{for}\ \ t\ge 0.
$$
Let $ \overline{u}_n = (4\pi/\alpha_0 - \delta)^{1/2} u_n / \|
u_n \|$ and using the inequality \re{inequa} with $s =
g(v_n)/C_\e$ and $t = \sqrt{\alpha_0}\, |\overline{u}_n|$,
{\allowdisplaybreaks
\begin{align*}
&(4\pi/\alpha_0 - \delta)^{1/2} \| u_n \| \le \int_{\R^2}g(v_n)|\overline{u}_n|\, \ud x \\
&=\frac{C_\e}{\sqrt{\alpha_0}}\int_{\R^2}\frac{g(v_n)}{C_\e}\sqrt{\alpha_0}|\overline{u}_n|\, \ud x \\
& \leq\frac{C_\e}{\sqrt{\alpha_0}}\int_{\{x \in\R^2 : g(v_n(x))/C_\e  \geq e^{1/\sqrt[3]{4}} \}}
\frac{g(v_n)}{C_\e} [\log(\frac{ g(v_n)}{C_\e})]^{1/2}\, \ud x  \\
&\ \ +\int_{\{x \in\R^2 :  g(v_n(x))/C_\e \leq e^{1/\sqrt[3]{4}} \}}g(v_n)|\overline{u}_n|\, \ud x+C_\e\sqrt{\alpha_0}\int_{\R^2}\overline{u}_n^2 ( e^{\alpha_0 \overline{u}_n^{2}} - 1 )\, \ud x\\
& \leq\sqrt{\frac{\alpha_0+\e}{\alpha_0}}\int_{\{x \in\R^2 : g(v_n(x))/C_\e  \geq e^{1/\sqrt[3]{4}} \}}
g(v_n)v_n\, \ud x+C_\e\sqrt{\alpha_0}\int_{\R^2}\overline{u}_n^2 ( e^{\alpha_0 \overline{u}_n^{2}} - 1 )\, \ud x   \\
&\ \ +\int_{\{x \in\R^2 :  g(v_n(x))/C_\e \leq e^{1/\sqrt[3]{4}} \}}g(v_n)|\overline{u}_n|\, \ud x\\
&\le \sqrt{\frac{\alpha_0+\e}{\alpha_0}}\int_{\R^2}g(v_n)v_n\, \ud x+I_{1,n}+I_{2,n}.
\end{align*}
}%
Recalling that $\overline{u}_n\rg0$ strongly in $L^s(\R^2)$ for any $s>2$. Since $\| \overline{u}_n \|^2 = 4\pi/\alpha_0 - \delta$, there exists $p>1$(close to 1) such that $p\al_0(4\pi/\alpha_0 - \delta)<4\pi$. Thus,
by the Pohozaev-Trudinger-Moser inequality, as $n\rg\iy$,
$$
I_{1,n}\le C_\e\sqrt{\al_0}\left(\int_{\R^2}|\overline{u}_n|^{2q}\right)^{1/q}\left(\int_{\R^2}( e^{p\alpha_0 \overline{u}_n^{2}}-1)\right)^{1/p}\rg0,
$$
where $1/p+1/q=1$, namely, $I_{1,n}=o_n(1)$. Note that by $(H1)$-$(H2)$, for any $\rho>0$, there exits $C_{\rho,\e}>0$ such that
$$
g(v_n(x))\le\rho|v_n(x)|+C_{\rho,\e}v_n^2,\ \ \mbox{for any}\  x\in\R^2\ \mbox{with}\ \ g(v_n(x))/C_\e \leq e^{1/\sqrt[3]{4}}.
$$
Then
\begin{align*}
I_{2,n}\le\int_{\R^2}(\rho|v_n|+C_{\rho,\e}v_n^2)|\overline{u}_n|\le\left[\rho\left(\int_{\R^2}|v_n|^2\right)^{1/2}+
C_{\rho,\e}\left(\int_{\R^2}|v_n|^4\right)^{1/2}\right]\left(\int_{\R^2}|\overline{u}_n|^2\right)^{1/2}.
\end{align*}
Recalling $v_n\rg0$ strongly in $L^4(\R^2)$,
$$
\limsup_{n\rg\iy}I_{2,n}\le C'\rho,
$$
where $C'>0$ is independent of $\rho$. By the arbitrary choice of $\rho$, $I_{2,n}=o_n(1)$. Hence,
\begin{equation}\label{barao4}
(4\pi/\alpha_0 - \delta)^{1/2} \| u_n \| \leq o_n(1) +
(1 + \frac{\e}{\alpha_0})^{1/2}\int_{\R^2 } g(v_n) v_n .
\end{equation}
Similarly, we have
\begin{equation}\label{CPV4}
(4\pi/\alpha_0 - \delta)^{1/2} \| v_n \| \leq o_n(1) +
(1 + \frac{\e}{\alpha_0})^{1/2} \int_{\R^2 } f(u_n) u_n.
\end{equation}
Note that
$$
\lan\Phi'(z_n),z_n\ran=2\int_{\R^2}(\nabla u_n \nabla v_n+V_0u_nv_n)-\int_{ \R^2 } f(u_n) u_n + \int_{ \R^2 } g(v_n) v_n=0
$$
and that by \re{y4.5} we get
$$\int_{ \R^2 } f(u_n) u_n+\int_{ \R^2 } g(v_n) v_n=2c_\ast+o_n(1).$$
It follows from (\ref{barao4})-(\ref{CPV4}) that
$$
(4\pi/\alpha_0 - \delta)^{1/2}(\| u_n \|_{H^1} + \| v_n \|_{H^1})\le 2c_\ast(1 + \frac{\e}{\alpha_0})^{1/2}+o_n(1).
$$
Since $c_\ast<4\pi/\alpha_0 - \delta$, for $\e > 0$ sufficiently small, as $n$ is large enough we have
$$
\| u_n \|_{H^1} + \| v_n \|_{H^1}\le 2(4\pi/\alpha_0 - \delta/2)^{1/2}.
$$

Then similarly as above, by the Trudinger-Moser inequality and $u_n\rg0$ strongly in $L^q(\R^2)$ for any $q>2$, we have $\int_{\R^2}g(v_n)u_n\rg0$, which implies by (\ref{unvn4}) that
$ u_n \rightarrow 0 $ strongly in $ H^1(\R^2)$. Thus, it follows from \re{y4.5} that $c_\ast=0$ and hence a contradiction and vanishing does not occur. 
\ep

\noindent As a consequence of Proposition \ref{nv}, up to a subsequence, there exist $\{y_n\}\subset\R^2$ and $z_0\not\equiv 0$ such that $z_n(\cdot+y_n)\rightharpoonup z_0$ in $E$ and $z_n(\cdot+y_n)\xrightarrow{a.e.}z_0$ in $\R^2$, as $n\rg\iy$.

\bo\lab{o11}
The weak limit $z_0$ is a critical point of $\Phi$.
\eo
\bp
By $(H1)$, there exist $a>0$ and $\al>\al_0$ such that
$$
|f(t)|\le a|t|+(e^{\al t^2}-1)\ \ \mbox{for all} \ \ t\in\R.
$$
Then by the Pohozaev-Trudinger-Moser inequality $f(\bar{u}_n)\in L_{loc}^1(\R^2)$ and $g(\bar{v}_n)\in L_{loc}^1(\R^2)$, where $\bar{z}_n=(\bar{u}_n,\bar{v}_n)=(u(\cdot+y_n),v(\cdot+y_n))$. From Lemma \ref{l3.3} and and Proposition \ref{o1} we get, as $n\rg\iy$
$$
\int_{\R^2}(f(\bar{u}_n)\varphi+g(\bar{v}_n)\phi)\rg\int_{\R^2}(f(u_0)\varphi+g(v_0)\phi),
$$
for any $(\varphi,\phi)\in C_0^\iy(\R^2)\times C_0^\iy(\R^2)$. Noting that $\Phi'(\bar{z}_n)=0$, it follows that
$$
\int_{\R^2}(\na u_0\na \phi+\na v_0\na \varphi+V_0u_0\phi+V_0v_0\varphi)\, \ud x=\int_{\R^2}(f(u_0)\varphi+g(v_0)\phi)\, \ud x,
$$
for any $(\varphi,\phi)\in C_0^\iy(\R^2)\times C_0^\iy(\R^2)$. Thus, $\Phi'(z_0)=0$ in $E$ and $z_0=(u_0,v_0)$ is a critical point of $\Phi$.
\ep
\vskip0.1in
\bo\lab{con}
$z_0\in \mathcal{S}$ and $z_n(\cdot+y_n)\longrightarrow z_0$ in $E$, as $n\rg\iy$, thus $\mathcal{S}$ is a compact set. 
\eo
\bp Thanks to the invariance of $\Phi$ by translation, let us write for simplicity $z_n$ in place of $z_n(\cdot+y_n)$ and let $z_n=(u_n,v_n)$, $z_0=(u_0,v_0)$. By $(H2)$, $f(s)s-2F(s)\ge0$ and $g(s)s-2G(s)\ge0$ for any $s\in\R$. Then by Fatou's Lemma,
{\allowdisplaybreaks
\begin{align}\lab{fatou}
c_\ast&=\Phi(z_n)-\frac{1}{2}\lan \Phi'(z_n),z_n\ran\nonumber\\
&=\lim_{n\rg\iy}\left(\int_{\R^2}\frac{1}{2}f(u_n)u_n-F(u_n)+\int_{\R^2}\frac{1}{2}g(u_n)u_n-G(u_n)\right)\nonumber\\
&\ge\limsup_{n\rg\iy}\int_{\R^2}\frac{1}{2}f(u_n)u_n-F(u_n)+\liminf_{n\rg\iy}\int_{\R^2}\frac{1}{2}g(u_n)u_n-G(u_n)\nonumber\\
&\ge\liminf_{n\rg\iy}\int_{\R^2}\frac{1}{2}f(u_n)u_n-F(u_n)+\liminf_{n\rg\iy}\int_{\R^2}\frac{1}{2}g(u_n)u_n-G(u_n)\\
&\ge\int_{\R^2}\frac{1}{2}f(u_0)u_0-F(u_0)+\int_{\R^2}\frac{1}{2}g(u_0)u_0-G(u_0)\nonumber\\
&=\Phi(z_0)-\frac{1}{2}\lan \Phi'(z_0),z_0\ran=\Phi(z_0).\nonumber
\end{align}}

On the other hand, since $z_0\not\equiv 0$ and $\Phi'(z_0)=0$ one has $\Phi(z_0)\ge c_\ast$. Therefore,
$z_0$ is a ground state solution of \eqref{q11}, namely, $z_0\in \mathcal{S}$.
\vskip0.1in
\noindent Next we prove that $z_n\rg z_0$ in $E$. By \re{fatou} and $\Phi(z_0)=c_\ast$ we have
\be\lab{ff}
\lim_{n\rg\iy}\int_{\R^2}\frac{1}{2}f(u_n)u_n-F(u_n)=\int_{\R^2}\frac{1}{2}f(u_0)u_0-F(u_0)
\ee
and
\be\lab{gg}
\lim_{n\rg\iy}\int_{\R^2}\frac{1}{2}g(v_n)v_n-G(v_n)=\int_{\R^2}\frac{1}{2}g(v_0)v_0-G(v_0).
\ee
By $(H2)$ we get 
$$
0\le\frac{\theta-2}{2}F(u_n)\le \frac{1}{2}f(u_n)u_n-F(u_n),\quad 0\le\frac{\theta-2}{2}G(v_n)\le \frac{1}{2}g(v_n)v_n-G(v_n), \quad n\geq 1
$$
and the Lebesgue dominated convergence theorem, together with \re{ff} and \re{gg} yields
\be\lab{fg1}
\lim_{n\rg\iy}\int_{\R^2}F(u_n)=\int_{\R^2}F(u_0),\quad  \lim_{n\rg\iy}\int_{\R^2}G(v_n)=\int_{\R^2}G(v_0).
\ee
Then, by \re{ff} and \re{gg} one has 
\be\lab{fg2}
\lim_{n\rg\iy}\int_{\R^2}f(u_n)u_n=\int_{\R^2}f(u_0)u_0,\quad \lim_{n\rg\iy}\int_{\R^2}g(v_n)v_n=\int_{\R^2}g(v_0)v_0.
\ee
Since $z_n,z_0\in S$, we have 
$$
\int_{\R^2}\na u_n\na v_n+V_0u_nv_n=c_\ast+\int_{\R^2}F(u_n)+G(v_n),
$$
$$
\int_{\R^2}\na u_0\na v_0+V_0u_0v_0=c_\ast+\int_{\R^2}F(u_0)+G(v_0).
$$
Thanks to \re{fg1},
\be\lab{fg3}
\lim_{n\rg\iy}\int_{\R^2}\na u_n\na v_n+V_0u_nv_n=\int_{\R^2}\na u_0\na v_0+V_0u_0v_0.
\ee
By $\lan\Phi'(u_n,v_n),(u_n,u_n)\ran=0$ and \re{fg2}-\re{fg3} we obtain 
\be\lab{fg4}
\int_{\R^2}|\na u_n|^2+V_0u_n^2=\int_{\R^2}f(u_0)u_0+g(v_n)u_n-\int_{\R^2}\na u_0\na v_0+V_0u_0v_0+o_n(1).
\ee
At the same time from $\lan\Phi'(u_n,v_n),(u_n,-u_n)\ran=0$ and $\lan\Phi'(u_0,v_0),(u_0,-u_0)\ran=0$, we have
\be\lab{fg5}
\int_{\R^2}f(u_n)u_n=\int_{\R^2}g(v_n)u_n,\quad \int_{\R^2}f(u_0)u_0=\int_{\R^2}g(v_0)u_0.
\ee
This implies by \re{fg2} that $\lim_{n\rg\iy}\int_{\R^2}g(v_n)u_n=\int_{\R^2}g(v_0)u_0$. As a consequence, by \re{fg4} we obtain 
$$
\lim_{n\rg\iy}\int_{\R^2}|\na u_n|^2+V_0u_n^2=\int_{\R^2}f(u_0)u_0+g(v_0)u_0-\int_{\R^2}\na u_0\na v_0+V_0u_0v_0.
$$
Recalling that $\lan\Phi'(u_0,v_0),(u_0,u_0)\ran=0$, namely 
$$\int_{\R^2}|\na u_0|^2+V_0u_0^2=\int_{\R^2}f(u_0)u_0+g(v_0)u_0-\int_{\R^2}\na u_0\na v_0+V_0u_0v_0,$$
which implies $$\lim_{n\rg\iy}\int_{\R^2}|\na u_n|^2+V_0u_n^2=\int_{\R^2}|\na u_0|^2+V_0u_0^2$$ and hence $u_n\rg u_0$ in $H^1(\R^2)$. Similarly, $v_n\rg v_0$ in $H^1(\R^2)$.
\ep

\noindent
Next we prove $(i), (iii)$ of Theorem \ref{Th2} through the following three steps:
\begin{itemize}
\item In Proposition \ref{bo1} we prove regularity, namely for any fixed $z=(u,v)\in \mathcal{S}$ we prove that $u,v\in L^{\iy}(\R^2)\cap C_{loc}^{1,\g}(\R^2)$ for some $\g\in(0,1)$;
\item In Proposition \ref{bo2} we prove that for any $\{z_n\}\subset \mathcal{S}$, $z_n=(u_n,v_n)$, for which there exists $y_n\in \mathbb{R}^2$ such that $z_n(\cdot +y_n)\to z_0\in \mathcal{S}$, one has 
$$\sup_{n\ge1}(\|u_n\|_\iy+\|v_n\|_\iy)<\iy;$$
\item Finally, in Proposition \ref{pro_apriori} we prove the following a priori estimates 
$$0<\inf_{z=(u,v)\in \mathcal{S}}\min\{\|u\|_\iy,\|v\|_\iy\}< \sup_{z=(u,v)\in \mathcal{S}}(\|u\|_\iy+\|v\|_\iy)<\iy.$$
\end{itemize}
\bo\lab{bo1}
Let $(u,v)\in \mathcal{S}$, then $u,v\in L^{\iy}(\R^2)\cap C_{loc}^{1,\g}(\R^2)$ for some $\g\in(0,1)$.
\eo
\bp For any $r>0$, let $B_1=B_r(0), B_2=B_{2r}(0)$. Noting that $u$ is a weak solution of the following problem
\be\lab{v0}
-\DD U+V_0U=g(v)\ \mbox{in}\ B_2,\ U-u\in H_0^1(B_2),
\ee
by the Pohozaev-Trudinger-Moser inequality one has $g(v)\in L^p(B_2)$ for all $p\ge2$. By the Calderon-Zygmund inequality, see e.g. \cite[Theorem 9.9]{GT}, one has $u\in W^{2,p}(B_2)$. It follows from classical interior $L^p$-estimates that
\be\lab{v1}
\|u\|_{W^{2,p}(B_1)}\le C\left(\|g(v)\|_{L^p(B_2)}+\|u\|_{L^p(B_2)}\right),
\ee
where $C$ only depends on $r,p$. Meanwhile, by the Sobolev embedding theorem, if $p>2$ we get that $u\in C^{1,\g}(\overline{B_1})$ for some $\g\in(0,1)$ and there exists $c$ (independent of $u$) such that
\be\lab{v2}\|u\|_{C^{1,\g}(\overline{B_1})}\le c\|u\|_{W^{2,p}(B_1)}.
\ee

\noindent Next we prove that $u$ vanishes at infinity, namely that for any $\dd>0$, there exists $R>0$ such that $|u(x)|\le \dd,\ \forall |x|\ge R$. Indeed, otherwise there exists $\{x_j\}\subset\R^2$ with $|x_j|\rg\iy$, as
$j\rg\iy$ and $\liminf_{j\rg\iy}|u(x_j)|>0$. Let $u_j(x)=u(x+x_j)$ and $v_j(x)=v(x+x_j)$, then $\|u_j\|=\|u\|$ and
\be\lab{v3}
-\DD u_j+V_0u_j=g(v_j),\quad  u_j\in H^1(\R^2).
\ee
We may assume $u_j \rightharpoonup u_0$ weakly in $H^1(\R^2)$, we claim that $u_0\not\equiv 0$. In fact, noting that $u_j$ is a weak solution of \eqref{v0} replacing $g(v)$ by $g(v_j)$, it follows from \eqref{v1} and
\eqref{v2} that, up to a subsequence, $u_j\rg u_0$ uniformly in $\overline{\om}$. Hence,
$$
u_0(0)=\liminf_{j\rg\iy}u_j(0)=\liminf_{j\rg\iy}u(x_j)\not=0,
$$
which implies that $u_0\not\equiv 0$. On the other hand, for any fixed $R>0$ and $j$ large enough, we have
\begin{align*}
\int_{\R^2}u^2 \, \ud x &\ge \int_{B_R(0)}u^2 \, \ud x +\int_{B_R(x_j)}u^2 \, \ud x\\
&=\int_{B_R(0)}u^2 \, \ud x+\int_{B_R(0)}u_j^2\, \ud x\\
&=\int_{B_R(0)}u^2\, \ud x +\int_{B_R(0)}u_0^2\, \ud x +o_j(1),
\end{align*}
where $o_j(1)\rg 0$, as $j\rg\iy$.
Since $R$ is arbitrary, we get $u_0\equiv 0$, which is a contradiction. Thus, $u(x)\rg 0$, as $|x|\rg\iy$. Moreover, since $u\in C(B_r)$ for any $r>0$, we have $u\in L^\iy(\R^2)$. Similarly, $v\in L^\iy(\R^2)$.
\ep
\bo\lab{bo2}
Let $z_n=(u_n,v_n)\subset \mathcal{S}$ such that  $\bar{z}_n=z_n(\cdot+y_n)\rg z_0=(u_0,v_0)\in \mathcal{S}$ in $E$, then
$$\sup_{n\ge1}(\|u_n\|_\iy+\|v_n\|_\iy)<\iy$$
\eo
\bp Let $\bar{u}_n=u(\cdot+y_n), \bar{v}_n=v_n(\cdot+y_n)$. Similarly as above, $\bar{u}_n$ is a weak solution of the following problem
\be\lab{vb0}
-\DD U+V_0U=g(\bar{v}_n)\ \mbox{in}\ B_2,\ U-\bar{u}_n\in H_0^1(B_2).
\ee Moreover, for any $p\ge2$ we have
\be\lab{vb1}
\|\bar{u}_n\|_{W^{2,p}(B_1)}\le C\left(\|g(\bar{v}_n)\|_{L^p(B_2)}+\|\bar{u}_n\|_{L^p(B_2)}\right),
\ee
where $C$ only depends on $r,p$. By the Sobolev embedding theorem, if $p>2$ we get $\bar{u}_n\in C^{1,\g}(\overline{B_1})$ for some $\g\in(0,1)$ and there exists $c$ (independent of $n$) such that
\be\lab{vb2}\|\bar{u}_n\|_{C^{1,\g}(\overline{B_1})}\le c\|\bar{u}_n\|_{W^{2,p}(B_1)}.
\ee
Then by \re{vb1}-\re{vb2}, we get
\be\lab{vb3}\|\bar{u}_n\|_{C^{1,\g}(\overline{B_1})}\le c\left(\|g(\bar{v}_n)\|_{L^p(\R^2)}+\|\bar{u}_n\|_{L^p(\R^2)}\right).
\ee

\noindent By $(H1)$, for $\beta>\al_0$ and some $C>0$, we have $|g(t)|\le C(|t|+\exp{(\beta t^2)}-1), t\in\R$. Recalling that $\bar{v}_n\rg v_0$ in
$H^1(\R^2)$, we next prove that
\be\lab{yf1}
\lim_{n\rg\iy}\int_{\R^2}|\exp(p\beta \bar{v}_n^2)-\exp(p\beta v_0^2)| \, \ud x=0.
\ee
In fact, since $v_0\in L^{\iy}(\R^2)$, there exists $c>0$ such that
\begin{align*}
&\int_{\R^2}|e^{(p\beta \bar{v}_n^2)}-e^{(p\beta v_0^2)}| \, \ud x \\
&\le c\int_{\R^2}e^{(2p\beta |\bar{v}_n-v_0|^2)}|\bar{v}_n^2-v_0^2| \ud x \\
&= c\int_{\R^2}[e^{(2p\beta |\bar{v}_n-v_0|^2)}-1]|\bar{v}_n^2-v_0^2| \, \ud x +o_n(1)\\
&\le c\left(\int_{\R^2}\left[e^{(4p\beta |\bar{v}_n-v_0|^2)}-1\right]\, \ud x \right)^{{1}/{2}}\left(\int_{\R^2}\left|\bar{v}_n^2-v_0^2\right|^2\, \ud x \right)^{{1}/{2}}+o_n(1),
\end{align*}
where $o_n(1)\rg 0$, as $n\rg\iy$. From $\|\bar{v}_n-v_0\|_1\rg0$, as $n\rg\iy$ and the Pohozaev-Trudinger-Moser inequality, there exists $C$ such that $$\int_{\R^2}\left[e^{(4p\beta
|\bar{v}_n-v_0|^2)}-1 \right]\, \ud x  \le C$$ as $n$ is large enough; thus \eqref{yf1} follows.

Recalling that $\bar{z}_n\rg z_0$ in $E$, by \re{yf1} $\|g(\bar{v}_n)\|_{L^p(\R^2)}\rg\|g(v_0)\|_{L^p(\R^2)}$ as $n\rg\iy$. Finally we have
\be\lab{vb4}
\sup_{n\ge1}\|\bar{u}_n\|_{C^{1,\g}(\overline{B_1})}<\iy.
\ee
\vskip0.1in

\noindent Next we prove that $\bar{u}_n(x)\rg0$, uniformly as $|x|\rg\iy$. It is enough to prove that for any $\dd>0$, there exists $R>0$ such that $|\bar{u}_n(x)|\le \dd,\ \forall n\ge1, |x|\ge R$. Suppose this does not occur, so that there exists
$\{x_n\}\subset\R^2$ with $|x_n|\rg\iy$, as $n\rg\iy$ and $\liminf_{n\rg\iy}|\bar{u}_n(x_n)|>0$. Let $\ti{u}_n(x)=\bar{u}_n(x+x_n)$ and $\ti{v}_n(x)=\bar{v}_n(x+x_n)$, then
\be\lab{vb5}
-\DD \ti{u}_n+V_0\ti{u}_n=g(\ti{v}_n),\quad  \ti{u}_n\in H^1(\R^2).
\ee
We may assume $\ti{u}_n \rightharpoonup \ti{u}_0$ weakly in $H^1(\R^2)$ and we claim $\ti{u}_0\not\equiv 0$. For any $n\ge1$, $\ti{u}_n$ is a weak solution to the following problem
\be\lab{vbb0}
-\DD U+V_0U=g(\ti{v}_n)\ \mbox{in}\ B_2,\ U-\ti{u}_n\in H_0^1(B_2).
\ee Moreover,
\be\lab{vbb1}
\|\ti{u}_n\|_{W^{2,4}(B_1)}\le C\left(\|g(\ti{v}_n)\|_{L^4(B_2)}+\|\ti{u}_n\|_{L^4(B_2)}\right)
\ee
where $C$ depends on $r$ only. At the same time, by the Sobolev embedding theorem, we get $\ti{u}_n\in C^{1,\g}(\overline{B_1})$ for some $\g\in(0,1)$ and there exists $c$ (independent of $n$) such that
\be\lab{vbb2}\|\ti{u}_n\|_{C^{1,\g}(\overline{B_1})}\le c\|\ti{u}_n\|_{W^{2,4}(B_1)}.
\ee
Then by \re{vbb1}-\re{vbb2}, we get
$$\|\ti{u}_n\|_{C^{1,\g}(\overline{B_1})}\le c\left(\|g(v_n)\|_{L^4(\R^2)}+\|u_n\|_{L^4(\R^2)}\right).
$$
Then similar to \re{vb4}, $\sup_{n\ge1}\|\ti{u}_n\|_{C^{1,\g}(\overline{B_1})}<\iy$. Hence up to a subsequence, $\ti{u}_n\rg \ti{u}_0$ uniformly in $\overline{B_1}$. Thus,
$$
\ti{u}_0(0)=\liminf_{n\rg\iy}\ti{u}_n(0)=\liminf_{n\rg\iy}u_n(x_n)\not=0,
$$
which implies that $\ti{u}_0\not\equiv 0$. On the other hand, for any fixed $R>0$ and $j$ large enough, we have
\begin{align*}
&o_n(1)+\int_{\R^2}u_0^2 \, \ud x=\int_{\R^2}\bar{u}_n^2 \, \ud x\\
&\ge \int_{B_R(0)}\bar{u}_n^2 \, \ud x +\int_{B_R(x_n)}\bar{u}_n^2 \, \ud x\\
&=\int_{B_R(0)}\bar{u}_n^2 \, \ud x+\int_{B_R(0)}\ti{u}_n^2\, \ud x\\
&=\int_{B_R(0)}u_0^2\, \ud x +\int_{B_R(0)}\ti{u}_0^2\, \ud x +o_n(1),
\end{align*}
where $o_n(1)\rg 0$, as $n\rg\iy$ and we have used the fact that $\bar{u}_n=u_n(\cdot+y_n)\rg u_0$ in $H^1(\R^2)$.
Since $R$ is arbitrary, we get $\ti{u}_0\equiv 0$, which is a contradiction. Thus, $\bar{u}_n(x)\rg 0$, uniformly as $|x|\rg\iy$, which immediately implies by \re{vb4} that
$\sup_{n\ge1}\|u_n\|_\iy=\sup_{n\ge1}\|\bar{u}_n\|_\iy<\iy$. Similarly, $\sup_{n\ge1}\|v_n\|_\iy<\iy$. 
\ep
\bo\label{pro_apriori} The following a priori estimates hold
\begin{equation}\label{aprioribound}
0<\inf_{z=(u,v)\in \mathcal{S}}\min\{\|u\|_\iy,\|v\|_\iy\}<\sup_{z=(u,v)\in \mathcal{S}}(\|u\|_\iy+\|v\|_\iy)<\iy
\end{equation}
\eo
\bp
The upper bound is a consequence of Proposition \ref{bo2} and the fact that $\mathcal{S}$ is compact. 

\noindent In order to prove the lower bound we argue by contradiction and thus assume $$\inf_{z=(u,v)\in \mathcal{S}}\min\{\|u\|_\iy,\|v\|_\iy\}=0.$$ 
Then,  there exists $\{z_n\}\subset \mathcal{S}$ such that, without loss of generality, $\|v_n\|_\iy\rg0$, as $n\rg\iy$. From 
$$
\int_{\R^2}|\na u_n|^2+V_0u_n^2=\int_{\R^2}g(v_n)u_n,
$$
by $(H1)$ we have
$$
\int_{\R^2}|\na u_n|^2+V_0u_n^2\leq o_n(1)\left(\int_{\R^2}v_n^2\right)^{1/2}\left(\int_{\R^2}u_n^2\right)^{1/2}
$$
and hence $u_n\rg0$ in $H^1(\R^2)$. From 
$$
\int_{\R^2}|\na v_n|^2+V_0v_n^2=\int_{\R^2}f(u_n)v_n\le\left(\int_{\R^2}v_n^2\right)^{1/2}\left(\int_{\R^2}[f(u_n)]^2\right)^{1/2},
$$
together with the fact $u_n\rg0$ in $H^1(\R^2)$ which implies $\int_{\R^2}[f(u_n)]^2\rg 0$, we have also $v_n\rg0$ in $H^1(\R^2)$. Finally, as $(u_n,v_n)\in \mathcal{S}$, we obtain a contradiction from the following 
$$
0<c_\ast=\lim_{n\rg\iy}\left(\int_{\R^2}\na u_n\na v_n+V_0u_nv_n-\int_{\R^2}F(u_n)+G(v_n)\right)=0
$$
\ep
In order to complete the proof of Theorem \ref{Th2} it remains to show that ground states vanish at infinity and that enjoy a suitable Pohozaev-type identity in the whole plane; we prove these results in Proposition \ref{vanishing_R} and \ref{stanislav} of the next Section.   

\subsection{Vanishing and Pohozaev-type identity}

\bo{\rm(Uniform vanishing)}\label{vanishing_R}
Let $x_z\in\R^2$ be a maximum point of $|u(x)|+|v(x)|$, $z=(u,v)\in \mathcal{S}$. Then $u(x+x_z)\to 0$ and $v(x+x_z)\rg 0$, as $|x|\rg\iy$, uniformly for any $(u,v)\in \mathcal{S}$.
\eo
\noindent In order to prove Proposition \ref{vanishing_R} we need the following technical lemma

\bl\lab{boo4}
For any $\{z_n\}\subset \mathcal{S}, z_n=(u_n,v_n)$, up to a subsequence, $z_n(\cdot+x_n)\rg z_1$ in $E$, as $n\rg\iy$, where $\{x_n\}\subset\R^2$ is such that  $|u_n(x_n)|+|v_n(x_n)|=\max_{x\in\R^2}(|u_n(x)|+|v_n(x)|).$
\el
\bp
We first claim that there exist $\mu>0$ and $R_1>0$ such that
\be\lab{nvv}
\lim_{n\rg\iy}\int_{B_{R_1}(x_n)}(u_n^2+v_n^2)\, \ud x\ge\mu.
\ee
Let us argue by contradiction, indeed if not, for some $\{z_n\}\subset\mathcal{S}$ and any $R>0$, we get
$$
\lim_{n\rg\iy}\int_{B_R(x_n)}(u_n^2+v_n^2)\, \ud x=0.
$$
Let $\hat{u}_n=u_n(\cdot+x_n)$ and $\hat{v}_n=v_n(\cdot+x_n)$, then $\hat{u}_n,\hat{v}_n\rg0$ in $L_{loc}^2(\R^2)$, as $n\rg\iy$. Similarly as above, $\hat{u}_n$ is a weak solution of the following problem
$$
-\DD U+V_0U=g(\hat{v}_n)\ \mbox{in}\ B_2,\ U-\hat{u}_n\in H_0^1(B_2).
$$ By standard elliptic regularity we get $\hat{u}_n\in C^{1,\g}(\overline{B_1})$ for some $\g\in(0,1)$ and there exists $c$ (independent of $n$) such that for $p>2$,
\be\lab{vbb3}\|\hat{u}_n\|_{C^{1,\g}(\overline{B_1})}\le c\left(\|g(\hat{v}_n)\|_{L^p(\R^2)}+\|\hat{u}_n\|_{L^p(\R^2)}\right).
\ee
By Proposition \ref{bo2}, $\bar{z}_n\rg z_0$ in $E$, by \re{yf1} $\|g(\hat{v}_n)\|_{L^p(\R^2)}=\|g(\bar{v}_n)\|_{L^p(\R^2)}\rg\|g(v_0)\|_{L^p(\R^2)}$, as $n\rg\iy$. Then we have
\be\lab{vbb4}
\sup_{n\ge1}\|\hat{u}_n\|_{C^{1,\g}(\overline{B_1})}<\iy,
\ee
which implies by $\hat{u}_n\rg0$ in $L^2(B_1)$ that $\hat{u}_n\rg0$ uniformly in $B_1$. In particular, $\hat{u}_n(0)=u_n(x_n)\rg0$. Similarly, we have $\hat{v}_n(0)=v_n(x_n)\rg0$. Finally we obtain
$$
\lim_{n\rg\iy}\max_{x\in\R^2}(|u_n(x)|+|v_n(x)|)=\lim_{n\rg\iy}(|u_n(x_n)|+|v_n(x_n)|)=0,
$$
which implies
$$
\lim_{n\rg\iy}\min\{\|u_n\|_\iy,\|v_n\|_\iy\}=0
$$
and thus a contradiction. 

Now by \re{nvv} $\lim_{n\rg\iy}\int_{B_{R_1}(0)}(\hat{u}_n^2+\hat{v}_n^2)\, \ud x\ge\mu$ which combined with the local compactness of the embedding $H^1(\R^2)\hookrightarrow L^2(\R^2)$, yields up to a subsequence, $z_n(\cdot+x_n)=(\hat{u}_n+\hat{v}_n)\rightharpoonup z_1\not={0}$ in $E$
and $z_n(\cdot+x_n)\to z_1$ a.e. in $\R^2$, as $n\rg\iy$. Then arguing as in Proposition \ref{o11}-\ref{con}, we get $z_1\in \mathcal{S}$ and $z_n(\cdot+x_n)\rg z_1$ in $E$, as $n\rg\iy$, and this completes the proof. 
\ep
\noindent{\it Proof of Proposition \ref{vanishing_R}.} 

Next let us prove that for any $\dd>0$, there exists $R>0$ such that $|u(x+x_z)|+|v(x+x_z)|\le \dd, |x|\ge R$ for any $z=(u,v)\in \mathcal{S}$, where $x_z\in\R^2$ is a maximum point of $|u(x)|+|v(x)|$. If not, there exist $z_n=(u_n,v_n)\in \mathcal{S}$ and $\{x_n\}\subset\R^2$ such that $|x_n|\rg\iy$ as $n\rg\iy$ and $$\liminf_{n\rg\iy}(|u_n(x_n+x_{z_n})|+|v_n(x_n+x_{z_n})|)>0,$$ where $x_{z_n}\in\R^2$ is a maximum point of $|u_n(x)|+|v_n(x)|$. Without loss of generality, we may assume $\liminf_{n\rg\iy}|u_n(x_n+x_{z_n})|>0$. Let $\ti{u}_n(x)=u_n(x+x_n+x_{z_n})$ and $\ti{v}_n(x)=v_n(x+x_n+x_{z_n})$. Assume $\ti{u}_n \rightharpoonup \ti{u}_0$ weakly in $H^1(\R^2)$, in the following we claim $\ti{u}_0\not\equiv 0$. Indeed, by Lemma \ref{boo4}, up to a subsequence, there exists $z\in\mathcal{S}$ such that $(u_n(\cdot+x_{z_n}),v_n(\cdot+x_{z_n}))\rg z$ strongly in $E$. Then as in the proof of the above Lemma, by the elliptic estimates, up to a subsequence, for some $\ti{u}_0\in H^1(\R^2)$ and $\g\in(0,1)$, $\ti{u}_n\rg \ti{u}_0$ in $C_{loc}^{1,\g}(\R^2)$, as $n\rg\iy$. Hence,
$$
\ti{u}_0(0)=\liminf_{n\rg\iy}\ti{u}_n(0)=\liminf_{n\rg\iy}u_n(x_n+x_{z_n})\not=0,
$$
which implies that $\ti{u}_0\not\equiv 0$. On the other hand, proceeding as in Proposition \ref{bo2}, we get $\ti{u}_0\equiv 0$, which is a contradiction.
\qed

\bo{\rm(Pohozaev-type identity)}\label{stanislav}
For any $z=(u,v)\in \mathcal{S}$, the following Pohozaev-type identity holds true
\begin{equation}\label{idpo}
\int_{\R^2}(F(u)+G(v)-V_0uv)\,\ud x=0.
\end{equation}
\eo
\bp
By the proof of Proposition \ref{bo1} we know $u,v\in W_{\loc}^{2,p}(\R^2)$ for any $p\ge2$. Then $\DD u=V_0u-g(v)$ a.e. in $\R^2$ and $\DD v=V_0v-f(u)$ a.e. in $\R^2$. Following\cite{Pucci, van} we get
\begin{align}\lab{poha}
&\oint_{\pl B_r}\na u\na v\cdot(x,{\bf n})\,\ud s-\oint_{\pl B_r}\left(\sum_{i,j=1}^2 x_j\left(\frac{\pl u}{\pl x_j}\frac{\pl v}{\pl x_i}+\frac{\pl v}{\pl x_j}\frac{\pl u}{\pl x_i}\right),{\bf n}\right)\,\ud s\\
&=2\int_{B_r}(V_0uv-F(u)-G(v))\,\ud x,\nonumber
\end{align}
where $B_r(0):=\{x\in\R^2:|x|<r\}, r>0$ and $\bf{n}$ is the outward normal of $\pl B_r$ at $x$. From $\na u,\na v\in L^2(\R^2)$, by virtue of the coarea formula, there exits $r_n$ such that $r_n\rg\iy$ and
$$
r_n\oint_{\pl B_{r_n}}\left|\frac{\pl u}{\pl x_j}\frac{\pl v}{\pl x_i}\right|\,\ud s\rg0,\,\,\hbox{for any}\,\,i,j=1,2.
$$
As a consequence as $n\rg\iy$,
$$
\left|\oint_{\pl B_{r_n}}\na u\na v\cdot(x,{\bf n})\,\ud s\right|\le r_n\oint_{\pl B_{r_n}}\left|\na u\na v\right|\,\ud s\rg0
$$
and hence 
$$
\oint_{\pl B_r}\left(\sum_{i,j=1}^2 x_j\left(\frac{\pl u}{\pl x_j}\frac{\pl v}{\pl x_i}+\frac{\pl v}{\pl x_j}\frac{\pl u}{\pl x_i}\right),{\bf n}\right)\,\ud s\rg0.
$$
Then, let $r=r_n$ in \re{poha} to get, as $n\rg\iy$, identity \eqref{idpo}.
\ep

\subsection{Sign and symmetry properties}\label{sign_s}
This section is devoted to proving Theorem \ref{sign}. To investigate positivity and radial symmetry of ground state solutions to \re{q11}, without loss
generality, throughout this section we assume that $f,g$ are odd symmetric functions.

\noindent Let $$\kappa:=\sup\{\|u\|_{\infty},\|v\|_{\infty}: (u,v)\in S\}<\iy$$ by Theorem \ref{Th2}. By $(H1)$ and $(H6)$, there exist small $a_0,b_0\in[0,1)$ and $k_1,k_2>0$ with $$k_1=\max_{a_0<|t|\le\kappa}|f(t)|/|t|^q,\,\, k_2=\max_{b_0<|t|\le\kappa}|g(t)|/|t|^p,$$ such that $f(t)\le t$, for $t\in[0,a_0]$ and $g(t)\le t$, for $t\in[0,b_0]$. Moreover, $f(a_0)=k_1a_0^q$ and $g(b_0)=k_2b_0^p$. In fact, if $\limsup_{t\rg0}|f(t)|/|t|^q<\iy$, we can choose $a_0=0$,  otherwise there exists $a_0\in(0,1)$ such that $f(a_0)/a_0^q=\max_{t\in[a_0,\kappa]}f(t)/t^q.$
Let
$$
f_k(t)=\left\{
\begin{array}{ll}
f(t),\ \ & \mbox{if}\ t\in[0,a_0]\\
\min\{f(t),k_1t^q\},\ \ \ & \mbox{if}\ t\in(a_0,\iy)
\end{array}
\right.
$$
and $f_k(t)=-f_k(-t)$ for $t\le0$ and similarly for $g$. Then, $f_k,g_k\in C(\R,\R)$ and $f_k(t)=f(t),g_k(t)=g(t)$ if $|t|\le\kappa$, $0<f_k(t)\le f(t),0<g_k(t)\le g(t)$ for all $t>0$. At the same time, there exists $\beta>0$ such that
\be\lab{gj}
\left\{
\begin{array}{ll}
|f_k(t)|\ge\beta|t|^q\text{ and }|g_k(t)|\ge\beta|t|^p,\ \ \mbox{for any}\ t\in\R\\
|f_k(t)|=|f(t)|\le|t|\,\,\,\mbox{if}\ |t|\le a_0,\ \ |g_k(t)|=|g(t)|\le|t|\,\ \mbox{if}\ |t|\le b_0\\
|f_k(t)|\le k_1|t|^q\,\,\mbox{if}\ |t|\ge a_0,\ \ |g_k(t)|\le k_2|t|^p\,\ \mbox{if}\ |t|\ge b_0.
\end{array}
\right.
\ee
Moreover, it is easy to check that $f_k, g_k$ satisfy $(H1)$, $(H4)$ and
\be\lab{am1}
0<2F_k(t)\le f_k(t)t,\,\,0<2G_k(t)\le g_k(t)t,\,\,t\not=0,
\ee
\be\lab{am2}
\lim_{|t|\rg\iy}\frac{F_k(t)}{t^2}=\iy,\,\,\lim_{|t|\rg\iy}\frac{G_k(t)}{t^2}=\iy,
\ee
where $F_k(t)=\int_0^tf_k(\tau)\,\ud\tau$ and $G_k(t)=\int_0^tg_k(\tau)\,\ud\tau$.

Now consider the truncated problem
\be\lab{qk1} \left\{
\begin{array}{ll}
&-\DD u+V_0u=g_k(v)\\ 
&-\DD v+V_0v=f_k(u)
\end{array}
\right. \ee
whose associated energy functional is
$$
\Phi_k(z):=\int_{\R^2}(\na u\na v+V_0uv)\,\ud x-\int_{\R^2}(F_k(u)+G_k(v))\,\ud x,\ \ z=(u,v)\in E.
$$
Recall the generalized Nehari Manifold
$$
\mathcal{N}_k:=\{z\in E\setminus E^-: \lan \Phi_k'(z),z\ran=0, \lan \Phi_k'(z),\vp\ran=0\ \mbox{for all}\ \ \vp\in E^-\}
$$
and the least energy
$$
c_\ast^k:=\inf_{z\in\mathcal{N}_k}\Phi_k(z).
$$
Noting that for any $(u,v)\in \mathcal{S}$, $(u,v)$ is a solution to \re{qk1}, hence $c_\ast^k\le c_\ast$. For $z\in E\setminus E^-$, set
$$
\hat{E}(z)=E^-\oplus\R^+z=E^-\oplus\R^+z^+.
$$
From \cite{DJJ,Szulkin,Weth} we have

\bl\lab{lk54.1}\
\begin{itemize}

\item [1)] For any $z\in \mathcal{N}_k$, $\Phi_k|_{\hat{E}(z)}$ has a unique maximum point which occurs exactly at $z$;

\item [2)] For any $z\in E\setminus E^-$, the set $\hat{E}(z)$ intersects $\mathcal{N}_k$ at exactly one point $\hat{m}_k(z)$, which is the unique global maximum point of $\Phi_k|_{\hat{E}(z)}$;

\item [3)]
$$
c_\ast^k:=\inf_{z\in E\setminus E^-}\max_{\omega\in\hat{E}(z)}\Phi_k(\omega).
$$
\end{itemize}
\el
\noindent From $0\le G_k(t)\le G(t)$ and $0\le F_k(t)\le F(t)$ for any $t\in\R$, we have
$$
c_\ast^k\ge\inf_{z\in E\setminus E^-}\max_{\omega\in\hat{E}(z)}\Phi(\omega)=c_\ast.
$$
thus $c_\ast^k=c_\ast>0$.

\noindent Next define 
$$
\hat{m}_k: z\in E\setminus E^-\mapsto\hat{m}_k(z)\in\hat{E}(z)\cap\mathcal{N}_k.
$$
There exists $\dd>0$ such that $\|z^+\|_\e\ge\dd$ for all $z\in\mathcal{N}_k$; in particular one has
$$
\|\hat{m}_k(z)^+\|_\e\ge\dd\ \ \ \mbox{for all}\ \ z\in E\setminus E^-.
$$
Moreover, for each compact subset $\mathcal{W}\subset E\setminus E^-$, there exists a constant $C_{\mathcal{W}}>0$ such that
$$
\|\hat{m}(z)\|\le C_{\mathcal{W}}\ \ \ \mbox{for all}\ \ z\in\mathcal{W}.
$$
\noindent Define
$$
S^+:=\{z\in E^+: \|z\|=1\},
$$
then, $S^+$ is a $C^1$-submanifold of $E^+$ and the tangent manifold of $S^+$ at $z\in S^+$ is given by
$$
T(S^+)=\{\omega\in E^+: (\omega,z)=0\}.
$$
Let
$$
m_k:=\hat{m}_k|_{S^+}: S^+\longrightarrow\mathcal{N}_k,
$$
then $\hat{m}_k$ is continuous and $m_k$ is a homeomorphism between $S^+$ and $\mathcal{N}_k$. Define
$$
\Psi_k: S^+\longrightarrow\R, \Psi_k(z):=\Phi_k(m_k(z)), z\in S^+
$$
then, by \cite[Corollary 4.3]{Weth} we have
\bo\lab{pk5.5}\noindent
\begin{itemize}

\item [1)] $\Psi_k\in C^1(S^+,\R)$ and
$$
\lan\Psi_k'(z),\omega\ran=\|m_k(z)^+\|\lan\Phi_k'(m_k(z)),\omega\ran,\ \ \mbox{for all}\ \ \omega\in T_z(S^+);
$$

\item [2)] If $\{\omega_n\}\subset S^+$ is a Palais-Smale sequence for $\Psi_k$, then $\{m_k(\omega_n)\}\subset \mathcal{N}_k$ is a Palais-Smale sequence for $\Phi_k$. Namely, if $\Psi_k(\omega_n)\rg d$ for some $d>0$ and $\|\Psi_k'(\omega_n)\|_\ast\rg 0$ as $n\rg\iy$, then $\Phi_k(m_k(\omega_n))\rg d$ and $\|\Phi_k'(m_k(\omega_n))\|\rg0$ as $n\rg\iy$, where
    $$
    \|\Psi_k'(\omega_n)\|_\ast=\sup_{\stackrel{\phi\in T_{\omega_n}(S^+)}{\|\phi\|=1}}\lan\Psi_k'(\omega_n),\phi\ran\ \ \mbox{and}\ \  \|\Phi_k'(m_k(\omega_n))\|=\sup_{\stackrel{\phi\in E}{\|\phi\|=1}}\lan\Phi_k'(m_k(\omega_n)),\phi\ran;
    $$

\item [3)] $\omega\in S^+$ is a critical point of $\Psi_k$ if and only if $m_k(\omega)\in \mathcal{N}_k$ is a critical point of  $\Phi_k$;

\item [4)] $\inf_{S^+}\Psi_k=\inf_{\mathcal{N}_k}\Phi_k$.
\end{itemize}
\eo
\noindent It follows from the Ekeland Variational Principle (see \cite[Theorem 3.1]{E}) that there exists $\{z_n^k\}\subset\mathcal{N}_k$ such that
\be\lab{pkss4}
\Phi_k(z_n^k)\rg c_\ast>0 \ \ \mbox{and}\ \ \Phi_k'(z_n^k)\rg 0,\ \ \mbox{as}\ \ n\rg\iy.
\ee
Next we prove that $\{z_n^k\}$ is uniformly bounded in $E$. Precisely we have the following
\bl\lab{o14}
There exists $C>0$ such that $\|z_n^k\|=\|(u_n^k,v_n^k)\|\le C$, for all $n\in\mathbb{N}$.
\el
\bp
Let $z_n^k=z_n^++z_n^-$, where $z_n^+\in E^+,\,\, z_n^-\in E^-$. Noting that $z_n^k\in\mathcal{N}_k$, we have $\|z_n^+\|^2\ge\|z_n^k\|^2/2$ for all $n\in\mathbb{N}$. Let $w_n^k=w_n^++w_n^-=z_n^k/\|z_n^k\|$, where $w_n^+\in\hat{E}(z_n^k)\subset E^+,w_n^-\in E^-$ and $w_n^+=(\ti{w}_n,\ti{w}_n)$, then $\|w_n^+\|^2\ge1/2$.  By Lemma \ref{lk54.1}, for some $R>2\sqrt{c_\ast}$, we have
\begin{align*}
c_\ast+o_n(1)&=\Phi_k(z_n^k)=\max_{w\in\hat{E}(z_n^k)}\Phi_k(w)\ge\Phi_k(R w_n^+)\\
&\ge R^2/4-\int_{\R^2}F_k(R\ti{w}_n)+G_k(R\ti{w}_n),
\end{align*}
which implies 
$$
\liminf_{n\rg\iy}\int_{\R^2}F_k(R\ti{w}_n)+G_k(R\ti{w}_n)>0.
$$
By Lions' Lemma, up to translations, $\ti{w}_n\rg w\not=0$ weakly in $H^1(R^2)$ as $n\rg\iy$. Assume that $w_n^k\rg (u,v)$ weakly in $H^1(R^2)$ as $n\rg\iy$, then $u+v=2w$. If $\|z_n^k\|\rg\iy$ as $n\rg\iy$, then $u_n^k(x)\rg\iy$ if $u(x)\not=0$ as $n\rg\iy$ and by Fatou's Lemma and \re{am2},
$$
\liminf_{n\rg\iy}\int_{\R^2}\left(\frac{F_k(u_n^k)}{\|z_n^k\|^2}+\frac{G_k(v_n^k)}{\|z_n^k\|^2}\right)=+\iy,
$$
which yields $\Phi_k(z_n^k)\rg-\iy$ as $n\rg\iy$. This is a contradiction and therefore $\{z_n^k\}$ stays bounded in $E$. 
\ep
Up to a subsequence, we may assume $z_n^k\rightharpoonup z^k$ weakly in $E$, as $n\rg\iy$. It is standard to check that $\Phi_k'(z^k)=0$.
\bo\lab{tk1}
The truncated problem \re{qk1} admits a ground state solution.
\eo
\bp
If $z^k\not=0$, then by \re{am2} and Fatou's Lemma one has 
\begin{align*}
c_\ast+o_n(1)&=\Phi_k(z_n^k)-\frac{1}{2}\lan \Phi_k'(z_n^k),z_n^k\ran\\
&=\int_{\R^2}\frac{1}{2}f_k(u_n^k)u_n^k-F_k(u_n^k)+\int_{\R^2}\frac{1}{2}g_k(u_n^k)u_n^k-G_k(u_n^k)\\
&\ge\int_{\R^2}\frac{1}{2}f_k(u^k)u^k-F_k(u^k)+\int_{\R^2}\frac{1}{2}g_k(u^k)u^k-G_k(u^k)+o_n(1)\\
&=\Phi_k(z^k)-\frac{1}{2}\lan \Phi_k'(z^k),z^k\ran+o_n(1)\\
&=\Phi_k(z^k)\ge c_\ast+o_n(1).
\end{align*}
from which $z^k$ is a ground state solution to \re{qk1}.

\noindent If $z^k=0$, we claim there exist $\nu>0$, $R_0>0$ and $\{y_n\}\subset\R^2$ such that
\be\lab{nonvanishing}
\lim_{n\rg\iy}\int_{B_{R_0}(y_n)}(|u_n^k|^2+|v_n^k|^2)\, \ud x\ge\nu.
\ee
Suppose the claim holds true and set $\ti{u}_n^k(\cdot):=u_n^k(\cdot+y_n)$ and $\ti{v}_n^k(\cdot):=v_n^k(\cdot+y_n)$, so that 
\be\lab{yk22}
\lim_{n\rg\iy}\int_{B_{R_0}(0)}(|\ti{u}_n^k|^2+|\ti{v}_n^k|^2)\, \ud x\ge\nu,
\ee
and $\Phi_k(\ti{z}_n^k)\rg c_\ast>0$ and $\Phi_k'(\ti{z}_n^k)\rg 0$, as $n\rg\iy$ where $\ti{z}_n^k=(\ti{u}_n^k,\ti{v}_n^k)$. Clearly $\{\ti{z}_n^k\}$ is bounded in $E$ and up to a subsequence, by \re{yk22} we may assume that $\ti{z}_n^k\rg \ti{z}^k\not=0$ weakly in $E$ to a ground state solution of \re{qk1}.

\noindent Hence let us prove by contradiction the claim \re{nonvanishing}. Indeed, if \re{nonvanishing} does not hold we have 
$$
\lim_{n\rg\iy}\sup_{y\in\R^2}\int_{B_R(y)}(|u_n^k|^2+|v_n^k|^2)\, \ud x=0\ \ \mbox{for all}\ \ R>0,
$$
then by Lions's Lemma, $u_n^k\rg0, v_n^k\rg0$ strongly in $L^s(\R^2)$ for any $s>2$. By $(H1)$ and \re{gj} we have
$$
\int_{\R^2}(|\na u_n^k|+V_0|u_n^k|)\,\ud x=\int_{\R^2}g_k(v_n^k)u_n^k\,\ud x\rg 0,\,\, n\rg\iy.
$$
Namely, $u_n^k\rg0$ strongly in $E$, as $n\rg\iy$. It follows that
$$
\int_{\R^2}(|\na v_n^k|+V_0|v_n^k|)\,\ud x=\int_{\R^2}f_k(u_n^k)v_n^k\,\ud x\rg 0,\,\, n\rg\iy.
$$
Namely, $v_n^k\rg0$ strongly in $E$, as $n\rg\iy$. So we get $c_\ast+o_n(1)=\Phi_k(z_n^k)\rg0$, as $n\rg\iy$, which is a contradiction.
\ep
\noindent Denote by $\mathcal{S}_k$ the set of of ground state solutions to system \re{qk1}, then $\mathcal{S}_k\not=\emptyset$. Similarly as above,
for any $z=(u,v)\in \mathcal{S}_k$, $u,v\in L^{\iy}(\R^2)\cap C_{loc}^{1,\g}(\R^2)$ for some $\g\in(0,1)$. Recalling that $c_\ast=c_\ast^k$,
we get $\mathcal{S}\subseteq \mathcal{S}_k$. In order to prove the reverse inclusion let us recall the following results from \cite{DJJ}

\bl\lab{l54.1}{\rm \cite{DJJ}} With the assumptions in Theorem \ref{Th2}, we have:
\begin{itemize}

\item [1)] for any $z\in \mathcal{N}$, $\Phi|_{\hat{E}(z)}$ admits a unique maximum point which is precisely at $z$;

\item [2)] for any $z\in E\setminus E^-$, the set $\hat{E}(z)$ intersects $\mathcal{N}$ at exactly one point $\hat{m}(z)$, which is the unique globally maximum point of $\Phi|_{\hat{E}(z)}$;

\item [3)]
$$c_\ast=\inf_{z\in E\setminus E^-}\max_{\omega\in\hat{E}(z)}\Phi(\omega).$$
\end{itemize}
\el
\noindent Let $m:=\hat{m}|_{S^+}: S^+\mapsto\mathcal{N}$ and
$$
\Psi: S^+\mapsto\R, \Psi(z):=\Phi(m(z)), z\in S^+,
$$
then $\hat{m}$ is continuous and $m$ is a homeomorphism between $S^+$ and $\mathcal{N}$. As in \cite{Weth}, $m$ is invertible and the inverse is given
by $$m^{-1}(z)=\frac{z^+}{\|z\|},\,\, z=z^++z^-\in\mathcal{N},\,\, z^+\in E^+,\,\ z^-\in E^-.$$
Similar to Proposition \ref{pk5.5}, we have
\bo\lab{pk5.6}\noindent
\begin{itemize}

\item [1)] $\Psi\in C^1(S^+,\R)$ and
$$
\lan\Psi'(z),\omega\ran=\|m(z)^+\|\lan\Phi'(m(z)),\omega\ran\ \ \mbox{for all}\ \ \omega\in T_z(S^+);
$$

\item [2)] If $\{\omega_n\}\subset S^+$ is a Palais-Smale sequence for $\Psi$, then $\{m(\omega_n)\}\subset \mathcal{N}$ is a Palais-Smale sequence for $\Phi$. Namely, if $\Psi(\omega_n)\rg d$ for some $d>0$ and $\|\Psi'(\omega_n)\|_\ast\rg 0$ as $n\rg\iy$, then $\Phi(m(\omega_n))\rg d$ and $\|\Phi'(m(\omega_n))\|\rg0$ as $n\rg\iy$, where
    $$
    \|\Psi'(\omega_n)\|_\ast=\sup_{\stackrel{\phi\in T_{\omega_n}(S^+)}{\|\phi\|=1}}\lan\Psi'(\omega_n),\phi\ran\ \ \mbox{and}\ \  \|\Phi'(m(\omega_n))\|=\sup_{\stackrel{\phi\in E}{\|\phi\|=1}}\lan\Phi'(m(\omega_n)),\phi\ran;
    $$

\item [3)] $\omega\in S^+$ is a critical point of $\Psi$ if and only if $m(\omega)\in \mathcal{N}$ is a critical point of  $\Phi$;

\item [4)] $\inf_{S^+}\Psi=\inf_{\mathcal{N}}\Phi$.
\end{itemize}
\eo

\bo\lab{sk}
$$\mathcal{S}_k=\mathcal{S}.$$
\eo
\bp
For any $z^k\in \mathcal{S}_k$, we know $z^k\in \mathcal{N}_k$, by Lemma \ref{lk54.1} $\Phi_k|_{\hat{E}(z)}$ admits a unique maximum point at $z^k$ and
$$
c_\ast^k:=\inf_{z\in E\setminus E^-}\max_{\omega\in\hat{E}(z)}\Phi_k(\omega)=\max_{\omega\in\hat{E}(z^k)}\Phi_k(\omega).
$$
Since $z^k\in E\setminus E^-$, by Lemma \ref{l54.1} the set $\hat{E}(z^k)$ intersects $\mathcal{N}$ just at one point $\hat{m}(z^k)$, which is the unique global maximum of $\Phi|_{\hat{E}(z^k)}$. Let $\hat{m}(z^k)=(\hat{u}^k,\hat{v}^k$), then by $0\le f_k(t)\le f(t)$ and $0\le g_k(t)\le g(t)$, for $t\ge0$ we have 
\begin{align*}
c_\ast^k&=\max_{\omega\in\hat{E}(z^k)}\Phi_k(\omega)\ge\Phi_k(\hat{m}(z^k))\\
&=\Phi(\hat{m}(z^k))+\int_{\R^2}[F(\hat{u}_k)-F_k(\hat{u}_k)]\,\ud x+\int_{\R^2}[G(\hat{v}_k)-G_k(\hat{v}_k)]\,\ud x\\
&=\max_{\omega\in\hat{E}(z^k)}\Phi(\omega)+\int_{\R^2}[F(\hat{u}_k)-F_k(\hat{u}_k)]\,\ud x+\int_{\R^2}[G(\hat{v}_k)-G_k(\hat{v}_k)]\,\ud x\\
&\ge\inf_{z\in E\setminus E^-}\max_{\omega\in\hat{E}(z)}\Phi(\omega)\ge c_\ast,
\end{align*}
which implies $F(\hat{u}_k(x))\equiv F_k(\hat{u}_k(x))$ and $G(\hat{v}_k(x))\equiv G_k(\hat{v}_k(x))$ for all $x\in\R^2$ and
$$\max_{\omega\in\hat{E}(z^k)}\Phi_k(\omega)=\Phi_k(\hat{m}(z^k))=\Phi(\hat{m}(z^k))=c_\ast.$$
Then $\Psi(m^{-1}(\hat{m}(z^k))):=\Phi(\hat{m}(z^k))=c_\ast$. Notice that $m^{-1}(\hat{m}(z^k))\in S^+$. Then, by Proposition \ref{pk5.6}, $m^{-1}(\hat{m}(z^k))$ is
a minimizer of $\Psi$ on the $C^1$-manifold $S^+$. Thus
$$
\lan\Psi'(m^{-1}(\hat{m}(z^k))),\omega\ran=0\ \ \mbox{for all}\ \ \omega\in T_{m^{-1}(\hat{m}(z^k))}(S^+).
$$
If follows from $3)$ of Proposition \ref{pk5.6} that $\Phi'(\hat{m}(z^k))=0$, which yields $\hat{m}(z^k)\in \mathcal{S}$. By uniqueness of the global maximum point
of $\Phi_k|_{\hat{E}(z^k)}$, we get $z^k=\hat{m}(z^k)$ and hence $z^k\in \mathcal{S}$. Therefore, $\mathcal{S}_k=\mathcal{S}$. 
\ep
In the last part of this section, in the spirit of \cite{dsr} we prove that $uv>0$ in $\R^2$ for any $z=(u,v)\in \mathcal{S}_k$.

\noindent Let $h(s):=g_k^{-1}(s)$ and $H$ denote the primitive function of $h$.
By \re{gj}, for some $c,C>0$,
\be\lab{gj2}
\left\{
\begin{array}{ll}
h(s)s\le C|s|^{(p+1)/p}\,\,\ &\mbox{for}\,\,s\in\R,\\
h(s)s\ge s^2/2\,\,\,&\mbox{if}\,\,|s|\le g(b_0),\\
h(s)s\ge c|s|^{(p+1)/p}\,\,\,&\mbox{if}\,\, |s|>g(b_0).
\end{array}
\right.
\ee
and clearly the same estimates hold for $H(s)$ as well. Consider the Schr\"odinger operator  $L:=-\DD +V_0$ and the Sobolev space $W^{2,(p+1)/p}(\R^2)$ endowed with the norm
$$\interleave{u\interleave}=\left(\int_{\R^2}|Lu|^{\frac{p+1}{p}}\,\ud x\right)^{\frac{p}{p+1}}.$$
The following embeddings hold
$$
W^{2,\frac{s+1}{s}}(\RN)\hookrightarrow L^r(\RN),\,\, \mbox{for any}\,\, r\ge\frac{s+1}{s},\,s>1,\,\,\mbox{if}\,\, s(N-2)\le 2,
$$
in particular $W^{2,(p+1)/p}(\R^2)\hookrightarrow L^2(\R^2)\cap L^{p+1}(\R^2)\cap L^{q+1}(\R^2)$. For $u\in W^{2,(p+1)/p}(\R^2)$, define
$$
J_k(u)=\int_{\R^2}H(Lu)-F_k(u)\,dx
$$
then $J_k$ is of class $C^1$ and
$$
\lan J_k'(u),\vp\ran=\int_{\R^2}(h(Lu)L(\vp)-f(u)\vp)\,\ud x,\,\, u,\vp\in W^{2,(p+1)/p}(\R^2).
$$
\bo\lab{equ}
$(u,v)\in E$ is a critical point of $\Phi_k$ if and only if $u$ is a critical point of $J_k$ and $v=h(Lu)$. Moreover, one has $\Phi_k(u,v)=J_k(u)$.
\eo
\noindent Define
$$
c_1(\R^2)=\inf_{u\in\mathcal{N}_J}J_k(u),\,\,\,\mbox{where}\,\,\,\mathcal{N}_J:=\{u\in W^{2,(p+1)/p}(\R^2)\setminus\{0\}: \lan J_k'(u),u\ran=0\},
$$
which under our assumptions might not be well defined. We overcome this difficulty by considering an approximation via bounded domains. Precisely, for any $R>0$ let us consider the problem
\be\lab{qk2} \left\{
\begin{array}{ll}
-\DD u+V_0u=g_k(v)\\
-\DD v+V_0v=f_k(u)
\end{array}
\right. \ee
$u,v\in H_0^1(B_R(0))$ whose associated energy functional is
$$
I_R(z):=\int_{B_R(0)}(\na u\na v+V_0uv)\,\ud x-(F_k(u)+G_k(v))\,\ud x,
$$
where $z=(u,v)\in E_R:=H_0^1(B_R(0))\times H_0^1(B_R(0))$.

\noindent We can define as above $E_R^+,E_R^-,\hat{E}_R(z)$ and
$$
\mathcal{N}_R:=\{z\in E_R\setminus E_R^-:\lan I_R'(z),z\ran=0,\,\,\lan I_R'(z),\phi\ran=0\,\,\, \mbox{for all}\,\,\phi\in E_R^-\}.
$$
Denote by $c_\ast(B_R(0))$ the corresponding least energy associated to the energy functional $I_R$. Similar to Lemma \ref{o14}, every Palais-Smale sequence for $I_R$ is bounded in $E_R$. Then $c_\ast(B_R(0))$ is the ground state critical level associated to $I_R$. Moreover,
$$
c_\ast(B_R(0))=\inf_{z\in E_R\setminus E_R^-}\max_{\omega\in\hat{E}_R(z)}I_R(\omega).
$$
\br
If $z=(u,v)\in\mathcal{N}_R$, we have $\lan I_R'(z),(\vp,-\vp)\ran=0$ for all $\vp\in H_0^1(B_R(0))$. In general, $\lan I_R'(z),(\vp,-\vp)\ran=0$ does not hold
for all $\vp\in H^1(\R^2)$. Then, $\mathcal{N}_R$ is not a subset of $\mathcal{N}$, so it is not clear if $c_\ast(B_R(0))$ is greater than $c_\ast$.
\er
\noindent Let $$X_R=W^{2,(p+1)/p}(B_R(0))\cap W_0^{1,(p+1)/p}(B_R(0))$$
endowed with the norm $$\interleave{u\interleave}=\left(\int_{B_R}|Lu|^{\frac{p+1}{p}}\,\ud x\right)^{\frac{p}{p+1}}$$ and
$$
J_R(u)=\int_{B_R(0)}H(Lu)-F_k(u)\,dx,\,\,\ u\in X_R.
$$

\bo
$z=(u,v)\in E_R$ is a critical point of $I_R$ if and only if $u$ is a critical point of $J_R$ and $v=h(Lu)$. Moreover, $I_R(u,v)=J_R(u)$.
\eo
\noindent Let
$$
\mathcal{N}_{J_R}:=\{u\in X_R\setminus\{0\}: \lan J_R'(u),u\ran=0\},\,\, c_1(B_R(0)):=\inf_{u\in\mathcal{N}_{J_R}}J_R(u).
$$
Notice that $\mathcal{N}_{J_R}$ might not be a $C^1-$manifold, so that we next borrow some ideas of \cite{Weth} to
overcome this difficulty and prove the existence of ground states corresponding to the functional $J_R$ on $\mathcal{N}_{J_R}$ for any $R$. Then by passing to
the limit, we show that $c_1(\R^2)$ is the ground state critical value.

\bl\lab{mountain}
For any $u\in X_R\setminus\{0\}$, $J_R(tu)\rg-\iy$, as $t\rg+\iy$ and the set $\R^+u$ intersects
$\mathcal{N}_{J_R}$ at exactly one point denoted by $\hat{m}_R(u)$, which is the unique global maximum point of $J_R(tu)$, for $t>0$. In particular,
$\hat{m}_R(u)=1$ if and only if $u\in\mathcal{N}_{J_R}$. Moreover, there exist $a_R,b_R>0$ such that
$$
\interleave{u\interleave}\ge a_R\,\,\mbox{for any}\,\,u\in\mathcal{N}_{J_R}\,\,\mbox{and}\,\,c_1(B_R(0))\ge b_R.
$$
\el
\bp
{\bf Step 1.} By \re{gj} and \re{gj2}, for any $u\in X_R\setminus\{0\}$ and $t>0$,
$$
J_R(tu)\le Ct^{(p+1)/p}\int_{B_R(0)}|Lu|^{(p+1)/p}-\frac{q+1}{q}\beta t^{q+1}\int_{B_R(0)}|u|^{q+1}\rg-\iy,\,\,t\rg+\iy,
$$
and for any $\g>0$ small, there exists $c_\g>0$ such that
\begin{align*}
J_R(tu)\ge& \frac{t^2}{2}\int_{\{|Lu|\le g(b_0)\}}|Lu|^2+ct^{(p+1)/p}\int_{\{|Lu|>g(b_0)\}}|Lu|^{(p+1)/p}\\
&-\g t^2\int_{B_R(0)}|u|^2-c_\g t^{q+1}\int_{B_R(0)}|u|^{q+1}>0,\,\,|t|\ll1,
\end{align*}
where
$$
\{|Lu|\le g(b_0)\}:=\{x\in B_R(0): |Lu(x)|\le g(b_0)\}.
$$
For any $u\in\mathcal{N}_{J_R}$, let $\theta(t)=J_R(tu)$, then $\theta(0)=0$ and $\theta'(1)=0$. Recalling that $g_k(s)/s$ is strictly increasing for $s>0$,
$h(s)/s$ is strictly decreasing for $s>0$. Obviously, $Lu=0$ if and only if $u=0$. Then for any $t>1$, thanks to $(H4)$, $(H6)$,
\begin{align*}
\theta'(t)&=\int_{B_R(0)}h(tLu)Lu-\int_{B_R(0)}f_k(tu)u\\
&=\int_{B_R(0)}h(t|Lu|)|Lu|-\int_{B_R(0)}f_k(t|u|)|u|\\
&=\int_{B_R(0)}\frac{h(t|Lu|)}{t|Lu|}t|Lu|^2-\int_{B_R(0)}\frac{f_k(t|u|)}{t|u|}t|u|^2\\
&<t\int_{B_R(0)}h(|Lu|)|Lu|-t\int_{B_R(0)}f_k(|u|)|u|\\
&=t\int_{B_R(0)}h(Lu)Lu-t\int_{B_R(0)}f_k(u)u=0.
\end{align*}
Similarly, $\theta'(t)>0$ for $t<1$. Namely, $J_R(u)=\max_{t\ge0}J_R(tu)$. Similarly, for any $u\in X_R\setminus\{0\}$, $J_R(tu)\rg-\iy$ as $t\rg+\iy$ and the set $\R^+u$ intersects $\mathcal{N}_{J_R}$ at exactly one point, which is the unique globally maximum point of $J_R(tu)$ for $t>0$.
\vskip0.1in
{\bf Step 2.} We prove that there exists $a_R>0$ such that
$$
\interleave{u\interleave}\ge a_R\,\,\mbox{for any}\,\,u\in\mathcal{N}_{J_R}.
$$
For any $u\in X_R\setminus\{0\}$, by \re{gj2} one has 
\begin{align}\lab{biao1}
\int_{B_R(0)}h(Lu)Lu&\ge\frac{1}{2}\int_{\{|Lu|\le g(b_0)\}}|Lu|^2+c\int_{\{|Lu|>g(b_0)\}}|Lu|^{(p+1)/p}\nonumber\\
&\ge\frac{1}{2}|B_R(0)|^{\frac{1-p}{1+p}}\left(\int_{\{|Lu|\le g(b_0)\}}|Lu|^{(p+1)/p}\right)^{2p/(p+1)}\\
&\ \ \ +c\int_{\{|Lu|>g(b_0)\}}|Lu|^{(p+1)/p}.\nonumber
\end{align}
Moreover, by $(H1)$, for any small $\g>0$, there exist $c_\g>0$ and $C>0$ (independent of $\g$) such that
\begin{align}\lab{biao2}
\int_{B_R(0)}f_k(u)u\le\int_{B_R(0)}\g u^2+c_\g |u|^{q+1}\le C\interleave{u\interleave}^2(\g+c_\g\interleave{u\interleave}^{q-1})
\end{align}
Here we used the embedding of $X_R$ into $L^r(B_R(0))$ for $r=2$ and $r=q+1$. By choosing
$$\g=2^{-\frac{4p+2}{p+1}}|B_R(0)|^{\frac{1-p}{1+p}}C^{-1},$$
and for any $u\in\mathcal{N}_{J_R}$, if $\interleave{u\interleave}^{q-1}\le\g c_\g^{-1}$, by \re{biao1} and \re{biao2},
\begin{align*}
&\frac{1}{4}|B_R(0)|^{\frac{1-p}{1+p}}\left(\int_{\{|Lu|\le g(b_0)\}}|Lu|^{(p+1)/p}\right)^{2p/(p+1)}+c\int_{\{|Lu|>g(b_0)\}}|Lu|^{(p+1)/p}\\
&\le C\g2^{2p/(p+1)}\left(\int_{\{|Lu|>g(b_0)\}}|Lu|^{(p+1)/p}\right)^{2p/(p+1)}.
\end{align*}
Since $u\not=0$, we have $\int_{\{|Lu|>g(b_0)\}}|Lu|^{(p+1)/p}>0$ and then
$$
\int_{\{|Lu|>g(b_0)\}}|Lu|^{(p+1)/p}\ge\left(\frac{c}{C\g2^{2p/(p+1)}}\right)^{\frac{p+1}{p-1}}>0.
$$
So that for any $u\in\mathcal{N}_{J_R}$ the following holds 
$$
\interleave{u\interleave}\ge\min\left\{(\g c_\g^{-1})^{\frac{1}{q-1}},\left(\frac{c}{C\g2^{2p/(p+1)}}\right)^{\frac{p}{p-1}}\right\}:=a_R>0.
$$
\vskip0.1in
{\bf Step 3.} We prove that there exists $b_R>0$ such that $c_1(B_R(0))\ge b_R.$ Obviously, $c_1(B_R(0))\ge0$. Assume by contradiction that there exists $\{u_n\}\subset\mathcal{N}_{J_R}$ such that $J_R(u_n)\rg0$, as $n\rg\iy$. We claim that $\{u_n\}$ is bounded in $X_R$. Indeed, if not we may assume $\interleave{u_n\interleave}\rg\iy$, as $n\rg\iy$. Let $v_n=u_n/\interleave{u_n\interleave}$ and assume that $v_n\rightharpoonup v$ weakly in $X_R$. If $v=0$, then by compactness of the embedding of $X_R$ into $L^r(B_R(0))$ for $r=2$ and $r=q+1$, we get $\int_{B_R(0)}F_k(v_n)\rg0$, as $n\rg\iy$. Then by Step 1,
\begin{align*}
J_R(u_n)=\max_{t\ge0}J_R(tu_n)\ge J_R(v_n)=\int_{B_R(0)}H(Lv_n)+o_n(1).
\end{align*}
Namely, $\int_{B_R(0)}H(Lv_n)=o_n(1)$. On the other hand, similar to \re{biao1},
\begin{align*}
\int_{B_R(0)}H(Lv_n)&\ge\frac{1}{2}|B_R(0)|^{\frac{1-p}{1+p}}\left(\int_{\{|Lv_n|\le g(b_0)\}}|Lv_n|^{(p+1)/p}\right)^{2p/(p+1)}\\
&\ \ \ +c\int_{\{|Lv_n|>g(b_0)\}}|Lv_n|^{(p+1)/p}.
\end{align*}
It follows that $v_n\rg0$ strongly in $X_R$, which contradicts the fact $\interleave{v_n\interleave}=1$. So $v\not=0$ and by \re{am2}, \re{gj2} and Fatou's Lemma,
$$
o_n(1)=\frac{J_R(u_n)}{\interleave{u_n\interleave}^{\frac{p+1}{p}}}\le C-\int_{B_R(0)}\frac{F_k(u_n)}{|u_n|^{(p+1)/p}}|v_n|^{(p+1)/p}\rg-\iy.
$$
This is a contradiction. Hence, $\{u_n\}$ is bounded in $X_R$. We may assume, up to a subsequence, $u_n\rightharpoonup u$ weakly in $X_R$ and strongly in $L^2(B_R(0))$. Noting that $h(t)/t$ is strictly decreasing for $t>0$, we have $0<h(t)t\le2H(t)$ for all $t\not=0$. Then by $(H2)$,
\begin{align*}
o_n(1)&=J_R(u_n)-\frac{1}{2}\lan J_R'(u_n),u_n\ran\\
&=\int_{B_R(0)}H(Lu_n)-\frac{1}{2}h(Lu_n)Lu_n+\frac{1}{2}\int_{B_R(0)}f_k(u_n)u_n-2F_k(u_n)\\
&\ge\frac{1}{2}\int_{B_R(0)}f_k(u_n)u_n-2F_k(u_n)\ge\frac{\theta-2}{2}\int_{\{x\in B_R(0): |u_n|\le a_0\}}F(u_n)\\
&\rg\frac{\theta-2}{2}\int_{\{x\in B_R(0): |u|\le a_0\}}F(u),\,\,\mbox{as}\,\,n\rg\iy.
\end{align*}
It follows that
$$
\int_{\{x\in B_R(0): |u|\le a_0\}}F(u)=0.
$$
Since $u\in X_R$, from elliptic regularity we get $u\in C^{0,2/(p+1)}(\overline{B_R(0)})$, which yields $u=0$. Analogously we get $\int_{B_R(0)}F_k(u_n)\rg0$, as $n\rg\iy$ and
\begin{align*}
\int_{B_R(0)}H(Lu_n)=J_R(u_n)+o_n(1)=o_n(1).
\end{align*}
Similar to \re{biao1},
\begin{align*}
\int_{B_R(0)}H(Lu_n)&\ge\frac{1}{2}|B_R(0)|^{\frac{1-p}{1+p}}\left(\int_{\{|Lu_n|\le g(b_0)\}}|Lu_n|^{(p+1)/p}\right)^{2p/(p+1)}\\
&\ \ \ +c\int_{\{|Lu_n|>g(b_0)\}}|Lu_n|^{(p+1)/p}.
\end{align*}
Thus $u_n\rg0$ strongly in $X_R$, which contradicts the fact $\interleave{u\interleave}\ge a_R$ for all $u\in\mathcal{N}_{J_R}$.
\ep
\noindent Define
$$
\hat{m}_R: u\in X_R\setminus\{0\}\mapsto\hat{m}_R(u)\in\R^+u\cap\mathcal{N}_{J_R}.
$$
Similar as in \cite{Szulkin}, we have the following 
\bl\lab{l2.4}
There exists $\dd>0$ such that $\interleave{u\interleave}\ge\dd$ for all $u\in\mathcal{N}_{J_R}$. In particular,
$$
\interleave{\hat{m}_R(u)\interleave}\ge\dd,\ \ \ \mbox{for all}\ \ u\in X_R\setminus\{0\}.
$$
Moreover, for each compact subset $\mathcal{W}\subset X_R\setminus\{0\}$, there exists a constant $C_{\mathcal{W}}>0$ such that
$$
\interleave{\hat{m}_R(u)\interleave}\le C_{\mathcal{W}},\ \ \ \mbox{for all}\ \ u\in\mathcal{W}.
$$
\el
\bp
By \re{gj2}, for any $u\in\mathcal{N}_{J_R}$, we have 
$$
b_1\le J_R(u)\le\int_{B_R(0)}H(Lu)\le C\interleave{u\interleave}^{p/(p+1)}.
$$
Thus, there exists $\dd>0$ such that $\interleave{u\interleave}\ge\dd$ for any $u\in\mathcal{N}_{J_R}$. Moreover, since $\hat{m}_R(u)=\hat{m}_R(u/\interleave{u\interleave})$ for any $u\not=0$, without loss generality,
we may assume $\mathcal{W}\subset S_R:=\{u\in X_R: \interleave{u\interleave}=1\}$. In the following, we claim that there exists $C_{\mathcal{W}}>0$ such that
\be\lab{claim}
\hbox{$J_R\le 0$ on $\R^+u\setminus B_{C_{\mathcal{W}}}(0)$, for all $u\in\mathcal{W}$,}
\ee
where $B_{C_{\mathcal{W}}}(0)=\{v\in X_R: \interleave{v\interleave}\le C_{\mathcal{W}}\}$.
If the claim \re{claim} is true, then noting that $J_R(\hat{m}_R(u))\ge b_1>0$ for all $0\not=u\in X_R$, we have
$\|\hat{m}_R(u)\|=\|\hat{m}_R(u/\interleave{u\interleave})\|\le C_{\mathcal{W}}$ for any $u\in \mathcal{W}$.
\vskip0.1in
\noindent So let us prove \re{claim}. Assume by contradiction that there exists $\{u_n\}\subset\mathcal{W}\subset S_R$ with $u_n\rg u$ strongly
in $\mathcal{W}$ and $\omega_n\in\R^+u_n$ with $\omega_n=t_nu_n$, $t_n\rg\iy$ such that $J_R(\omega_n)\ge0$, as $n\rg\iy$.
For $n$ large enough, by \re{gj2} one has 
\begin{align}\lab{y3}
0\le\frac{J_R(\omega_n)}{\interleave{\omega_n\interleave}^{(p+1)/p}}&\le C
-\int_{B_R(0)}\frac{F_k(t_nu_n)}{|t_nu_n|^{(p+1)/p}}|u_n|^{(p+1)/p}.
\end{align}
Noting that $u_n\xrightarrow{a.e.}u\not=0$, it follows from Fatou's Lemma and \re{y3} that  $$\frac{J_R(\omega_n)}{\interleave{\omega_n\interleave}^{(p+1)/p}}\rg-\iy$$
 as $n\rg\iy$, which is a contradiction. 
\ep
\noindent Let $m_R:=\hat{m}_R|_{S_R}: S_R\longrightarrow\mathcal{N}_{J_R}$ and
$$
K: S_R\longrightarrow \R,\quad  K(u):=J_R(m_R(u)), u\in S_R,
$$
then $\hat{m}_R$ is continuous and $m_R$ is a homeomorphism between $S_R$ and $\mathcal{N}_{J_R}$. 
\bo\lab{pkr5.6}\noindent
\begin{itemize}

\item [1)] $K\in C^1(S_R,\R)$ and
$
\lan K'(u),\omega\ran=\|m_R(u)\|\lan J_R'(m_R(u)),\omega\ran$, for all $\omega\in T_u(S_R)$;

\item [2)] If $\{\omega_n\}\subset S_R$ is a Palais-Smale sequence for $K$, then $\{m_R(\omega_n)\}\subset \mathcal{N}_{J_R}$
 is a Palais-Smale sequence for $J_R$. Namely, if $K(\omega_n)\rg d$ for some $d>0$ and $\|K'(\omega_n)\|_\ast\rg 0$, as $n\rg\iy$,
 then $J_R(m_R(\omega_n))\rg d$ and $\|J_R'(m_R(\omega_n))\|\rg0$, as $n\rg\iy$, where
    $$
    \|K'(\omega_n)\|_\ast=\sup_{\stackrel{\phi\in T_{\omega_n}(S_R)}{\interleave{\phi\interleave}=1}}\lan K'(\omega_n),\phi\ran\ \
    \mbox{and}\ \  \|J_R'(m_R(\omega_n))\|=\sup_{\stackrel{\phi\in X_R}{\interleave{\phi\interleave}=1}}\lan J_R'(m_R(\omega_n)),\phi\ran.
    $$

\item [3)] $\omega\in S_R$ is a critical point of $K$ if and only if $m_R(\omega)\in \mathcal{N}_{J_R}$ is a critical point of  $J_R$;

\item [4)] $\inf_{S_R}K=\inf_{\mathcal{N}_{J_R}}J_R$.
\end{itemize}
\eo
\bl
For any $R>0$, $c_1(B_R(0))\ge c_\ast(B_R(0)).$
\el
\bp Observing that $S_R$ is a $C^1$-manifold in $X_R$, by virtue of the Ekeland variational principle (see \cite[Theorem 3.1]{E}), there exists $\{u_n\}\subset\mathcal{N}_{J_R}$
such that
\be\lab{pkrss4}
J_R(u_n)\rg c_1(B_R(0))>0 \ \ \mbox{and}\ \ J_R'(u_n)\rg 0,\ \ \mbox{as}\ \ n\rg\iy.
\ee
It is standard to show that $\{u_n\}$ is bounded in $X_R$, thus up to a subsequence, $u_n\rg u$ weakly in $X_R$, as $n\rg\iy$. By means of the compactness
of $X_R\hookrightarrow L^{r}(B_R(0))$ for any $r\ge(p+1)/p$, $u_n\rg u$ strongly in $L^{q+1}(B_R(0))$. Then
\begin{align}\lab{nonzero}
\liminf_{n\rg\iy}\int_{B_R(0)}h(Lu_n)Lu_n=\liminf_{n\rg\iy}\int_{B_R(0)}f(u_n)u_n=\int_{B_R(0)}f(u)u.
\end{align}
By \re{gj2}, we also have
$$
\int_{B_R(0)}h(Lu_n)Lu_n\ge\frac{1}{2}\int_{|Lu_n|\le g(b_0)}|Lu_n|^2+c\int_{|Lu_n|>g(b_0)}|Lu_n|^{(p+1)/p}.
$$
We claim that $u\not\equiv0$. Indeed, otherwise by \re{nonzero} we get 
$$
\lim_{n\rg\iy}\int_{|Lu_n|\le g(b_0)}|Lu_n|^2=0\,\,\mbox{and}\,\,\ \lim_{n\rg\iy}\int_{|Lu_n|>g(b_0)}|Lu_n|^{(p+1)/p}=0.
$$
Hence 
\begin{align*}
\lim_{n\rg\iy}&\int_{B_R(0)}|Lu_n|^{(p+1)/p}\le\lim_{n\rg\iy}\int_{|Lu_n|>g(b_0)}|Lu_n|^{(p+1)/p}\\
&+\lim_{n\rg\iy}\left(\int_{|Lu_n|\le g(b_0)}|Lu_n|^2\right)^{(p+1)/(2p)}|B_R(0)|^{(p-1)/(2p)}
\rg0&
\end{align*}
as $n\to\infty$, which implies $J_R(u_n)\rg0$, as $n\rg\iy$. This is a contradiction.

\noindent Next let $u_0=\hat{m}_R(u)u$ and $v_n=\hat{m}_R(u)u_n$. By $(H7)$, $H$ is convex. Therefore
$$
\liminf_{n\rg\iy}\int_{B_R(0)}H(Lv_n)\ge\int_{B_R(0)}H(Lu_0)\,\,\mbox{and}\,\, \lim_{n\rg\iy}\int_{B_R(0)}F(v_n)=\int_{B_R(0)}F(u_0).
$$
As $u_0\in\mathcal{N}_{J_R}$ one the one hand on has 
\begin{align*}
\liminf_{n\rg\iy}J_R(v_n)\ge\int_{B_R(0)}H(Lu_0)-\int_{B_R(0)}F(u_0)\ge c_1(B_R(0)).
\end{align*}
On the other hand, it follows from Lemma \ref{mountain} and $u_n\in\mathcal{N}_{J_R}$ the following 
$$
\liminf_{n\rg\iy}J_R(v_n)\le\liminf_{n\rg\iy}\max_{t\ge0}J_R(tu_n)=\liminf_{n\rg\iy}J_R(u_n)=c_1(B_R(0)).
$$
and in turn $J_R(u_0)=c_1(B_R(0))$. By Proposition \ref{pkrss4} $J_R'(u_0)=0$ and by Proposition \ref{equ},
$(u_0,v_0)$ is a nontrivial critical point of $I_R$, namely $(u_0,v_0)\in\mathcal{N}_R$ where $v_0=h(Lu_0)$. Finally,
$$
c_\ast(B_R(0))\le I_R(u_0,v_0)=J_R(u_0)=c_1(B_R(0)).
$$
\ep
Similar as in \cite{dsr}, one can prove the reversed inequality to get the following 
\bl\lab{crr}
For any $R>0$,
$$c_\ast(B_R(0))=c_1(B_R(0)).$$
\el
\bl\lab{signr}
Let $(u_R,v_R)$ be any ground state for the functional $I_R$, then $u_Rv_R>0$ in $B_R(0)$.
\el
\bp
Recalling that $\mathcal{S}=\mathcal{S}_k$, it is enough to
prove $uv>0$ in $\R^2$ for any $(u,v)\in\mathcal{S}_k$. For any $R>0$ and any ground state $(u_R,v_R)$ for the functional $I_R$,
by Lemma \ref{crr} and Proposition \ref{equ}, $u_R$ is a ground state for the functional $J_R$. Let $\omega=L^{-1}(|Lu_R|)$, then $\omega>0$
and $\omega\ge|u_R|$. Moreover, $\lan J_R'(t\omega),\omega\ran=0$, where $t=\hat{m}_R(\omega)>0$. On the other hand,
\begin{align*}
c_1(B_R(0))&\le J_R(t\omega)=J_R(tu_R)+\int_{B_R(0)}F_k(t|u_R|)-F_k(t\omega)\\
&\le c_1(B_R(0))+\int_{B_R(0)}F_k(t|u_R|)-F_k(t\omega).
\end{align*}
So that
$\int_{B_R(0)}F_k(t|u_R|)-F_k(t\omega)\ge0$. It follows from $(H7)$ that $|u_R|=\omega>0$. If $u_R>0$ in $B_R(0)$, then by means of the maximum principle, $v_R>0$
in $B_R(0)$ and $u_Rv_R>0$ in $B_R(0)$. Similarly, if $u_R<0$ in $B_R(0)$, $u_Rv_R>0$ in $B_R(0)$.
\ep
\noindent As a consequence of Lemma \ref{crr} and Lemma \ref{signr}, see also \cite[Remark 4.11]{dsr}, we have
\bl\lab{compare}
The map $R\mapsto c_\ast(B_R(0))$ is decreasing for $R>0$.
\el
\bl\lab{cast}
For any $R>0$, we have $c_\ast(B_R(0))\ge c_\ast(\R^2)$.
\el
\bp
For any $R>0$, let $z_R=(u_R,v_R)$ be a ground state solution of $I_R$. Namely, $I_R(z_R)=c_\ast(B_R(0))$ and $I_R'(z_R)=0$.
We extend $z_R\in E_R$ to $z_R\in E$ by zero extension outside $B_R(0)$. Then, as in Lemma \ref{o14}, $\{z_R\}$ turns out to be 
bounded in $E$. Up to a subsequence, we may assume $z_R\rightharpoonup z_0$ weakly in $E$, as $R\rg\iy$, then $z_0=(u_0,v_0)\in E$ is a nonnegative
solution to \re{qk1}, namely $\Phi_k'(z_0)=0$.

If $z_0\not=0$, by $(H2)$ and Fatou's Lemma, we have for any $r\le R$,
\begin{align*}
c_\ast(B_r(0))&\ge\lim_{R\rg\iy}c_\ast(B_R(0))=\lim_{R\rg\iy}\left(I_R(z_R)-\frac{1}{2}\lan I_R'(z_R),z_R\ran\right)\\
&=\lim_{R\rg\iy}\left(\int_{B_R(0)}\frac{1}{2}f_k(u_R)u_R-F_k(u_R)+\int_{B_R(0)}\frac{1}{2}g_k(v_R)v_R-G_k(v_R)\right)\\
&\ge\int_{\R^2}\frac{1}{2}f_k(u_0)u_0-F_k(u_0)+\int_{\R^2}\frac{1}{2}g_k(v_0)v_0-G_k(v_0)\\
&=\Phi_k(z_0)-\frac{1}{2}\lan \Phi_k'(z_0),z_0\ran=\Phi_k(z_0)\ge c_\ast(\R^2).
\end{align*}
If $z_0=0$, then $\{z_R\}$ satisfies one of the following alternatives:
\begin{itemize}

\item [(1)] ({\it Vanishing})
$$
\lim_{R\rg\iy}\sup_{y\in\R^2}\int_{B_r(y)}(u_R^2+v_R^2)\, \ud x=0,\ \ \mbox{for all}\ \ r>0;
$$

\item [(2)] ({\it Nonvanishing}) there exist $\nu>0$, $r_0>0$ and $\{y_R\}\subset\R^2$ such that
$$
\lim_{R\rg\iy}\int_{B_{r_0}(y_R)}(u_R^2+v_R^2)\, \ud x\ge\nu.
$$
\end{itemize}
As in Proposition \ref{nv}{\it Vanishing} does not occur. So let $\ti{u}_R:=u_R(\cdot+y_R)$ and $\ti{v}_R:=v_R(\cdot+y_R)$,
then $\ti{z}_R=(\ti{u}_R,\ti{v}_R)$ is bounded in $H^1(\R^2)$ and $\ti{z}_R\rightharpoonup \ti{z}_0\not=0$ weakly in $H^1(\R^2)$.
Moreover, let $\ti{z}_0=(\ti{u}_0,\ti{v}_0)$, we know $\ti{u}_0,\ti{v}_0$ are nonnegative. Obviously, $|y_R|\le R+r_0$. Assume that, up to a rotation,
$y_R/|y_R|\rg (0,-1)\in\R^2$ and $(\ti{u}_0,\ti{v}_0)\in H_0^1(\Omega)\times H_0^1(\Omega)$ satisfies
\be\lab{qkj2} \begin{cases}
-\DD \ti{u}_0+V_0\ti{u}_0=g_k(\ti{v}_0)\\
-\DD \ti{v}_0+V_0\ti{v}_0=f_k(\ti{u}_0)
\end{cases}
\ee
where $\Omega=\R^2$ or $\Omega=\{(x_1,x_2)\in\R^2:x_2>d\}$, where $d:=\liminf_{R\rg\iy}\mbox{dist}(y_R,\pa B_R(0))$. If $\Omega=\R^2$ the proof follows. If $\Omega=\{(x_1,x_2)\in\R^2:x_2>d\}$, then by the Hopf Lemma,
$\pa\ti{u}_0/\pa \eta<0$ and $\pa\ti{v}_0/\pa \eta<0$ on $\pa\Omega$, where $\eta$ is the outward pointing unit normal to $\pa\Omega$. Finally from the Pohozaev type identity proved in \cite[Proposition 1.2]{Lions}(see also \cite[Lemma 3.1]{Pisto}) one actually has 
$$
\int_{\pa\Omega}\frac{\pa\ti{u}_0}{\pa n}\frac{\pa\ti{v}_0}{\pa n}=0,
$$
which is a contradiction.
\ep

\bp[Proof of Theorem \ref{sign} completed]
Thanks to Lemma \ref{cast} any ground state solution $(u,v)$ to \re{q11} does not change sign.
Assume $u>0$ and $v>0$ in $\R^2$. Setting 
$$
f_1(u,v)=g(v)-V_0u\,\,\,\mbox{and}\,\,\ f_2(u,v)=f(u)-V_0v,
$$
as a consequence of \cite[Theorem 1]{Sirakov1} and $(H1)$, $(u,v)$ is radially symmetric and strictly decreasing with respect to the same point, which we denote by $x_0$. Clearly, $\DD u(x_0)\le0$ and $\DD v(x_0)\le0$. To complete the proof of Theorem \ref{sign}, we next prove that actually $\DD u(x_0)<0$ and $\DD v(x_0)<0$. Indeed, if not, without loss of generality we may assume $\DD u(x_0)=0$ and then $g(v(x_0))=V_0u(x_0)$. Let $u_1=u-u(x_0)$, then $u_1(x)\le0$ in $\R^2$ and
\begin{align*}
&-\DD u_1=-\DD u=g(v)-V_0u\\
&\le g(v(x_0))-V_0u(x_0)-V_0u_1\\
&=-V_0u_1.
\end{align*}
Namely, $-\DD u_1+V_0 u_1\le 0$ in $\R^2$. Noting that $u_1(0)=0$, by the maximum principle, $u_1\equiv0$ in $\R^2$, which is a contradiction. Therefore, $\DD u(x_0)<0$. Similarly, one has $\DD v(x_0)<0$ as well. {Finally, by Proposition \ref{vanishing_R}, $u(x+x_z), v(x+x_z)\rg 0$, as $|x|\rg\iy$ uniformly for any $z=(u,v)\in \mathcal{S}$. Since $u,v$ do not change the sign, using the maximum principle, we conclude that there exist $C,c>0$, independent of $z=(u,v)\in \mathcal{S}$, such that $$|D^{\al}u(x)|+|D^{\al}v(x)|\le C\exp(-c|x-x_0|),\,x\in \R^2,\,|\al|=0,1$$ }
\ep

\s{Proof of Theorem \ref{Th1}}\label{semiclassical_s}
\renewcommand{\theequation}{5.\arabic{equation}}

\subsection{Functional setting} By setting $u(x)=\vp(\e x),v(x)=\psi(\e x)$ and $V_\e(x)=V(\e x)$, \re{q1} is equivalent to
\be\lab{q51} \left\{
\begin{array}{ll}
-\DD u+V_\e(x)u=g(v)\\
-\DD v+V_\e(x)v=f(u)
\end{array}
\right. \ee
We next consider \re{q51}. Let $H_\e$ be the completion of $C_0^\iy(\R^2)$ with respect to the inner product 
$$
(u,v)_{1,\e}:=\int_{\R^2}\na u\na v+V_\e(x)uv$$
and the norm 
$$\|u\|_{1,\e}^2:=(u,u)_{1,\e}, u,v\in H_\e.
$$
Let $E_\e:=H_\e\times H_\e$ with the inner product
$$
(z_1,z_2)_\e:=(u_1,u_2)_{1,\e}+(v_1,v_2)_{1,\e},\ \ z_i=(u_i,v_i)\in E_\e,\: i=1,2.
$$
and the norm $\|z\|_\e^2=\|(u,v)\|_\e^2=\|u\|_{1,\e}^2+\|v\|_{1,\e}^2.$ We have the orthogonal space decomposition $E_\e=E_\e^+\oplus E_\e^-$, where
$$
E_\e^+:=\{(u,u)\,|\, u\in H_\e\}\ \ \ \mbox{and}\ \ \ E_\e^-:=\{(u,-u)\,|\,u\in  H_\e\}.
$$
For each $z=(u,v)\in E_\e$, $$z=z^++z^-=((u+v)/2,(u+v)/2)+((u-v)/2,(v-u)/2).$$
Weak solutions of \re{q51} are critical points of the associated energy functional
$$
\Phi_\e(z):=\int_{\R^2}\na u\na v+V_\e(x)uv-I(z),\ \ z=(u,v)\in E_\e,
$$
where $I(z)=\int_{\R^2}F(u)+G(v)$. Then $\Phi_\e\in C^1(E,\R)$
and
$$
\lan\Phi_\e'(z),w\ran=\int_{\R^2}(\na u\na w_2+\na v\na w_1+V_\e(x)uw_2+V_\e(x)vw_1)-\int_{\R^2}(f(u)w_1+g(v)w_2),
$$
for all $z=(u,v),w=(w_1,w_2)\in E_\e$. Moreover, $\Phi_\e$ can be rewritten as follows
\be\lab{y51}
\Phi_\e(z):=\frac{1}{2}\|z^+\|_\e^2-\frac{1}{2}\|z^-\|_\e^2-I(z).
\ee
We know that if $z\in E_\e$ is a nontrivial critical point of $\Phi_\e$, then $z\in E_\e\setminus E_\e^-$. In the spirit of \cite{Szulkin}, we define the generalized Nehari Manifold
$$
\mathcal{N}_\e:=\{z\in E_\e\setminus E_\e^-: \lan \Phi_\e'(z),z\ran_\e=0, \lan \Phi_\e'(z),\vp\ran_\e=0\ \mbox{for all}\ \ \vp\in E_\e^-\}.
$$
Let
$$
c_\e:=\inf_{z\in\mathcal{N}_\e}\Phi_\e(z),
$$
then $c_\e$ is the least energy for system \re{q51}, the so-called ground state level.

\noindent For $z\in E_\e\setminus E_\e^-$, set
$$
\hat{E}_\e(z)=E_\e^-\oplus\R^+z=E_\e^-\oplus\R^+z^+,
$$
where $\R^+z:=\{tz: t\ge0\}$. From \cite{DJJ,Szulkin,Weth} we have the following properties of $\mathcal{N}_\e$, which will be used later.

\bl\lab{l5.1} Under the assumptions in Theorem \ref{Th1}, we have:
\begin{itemize}

\item [1)] for any $z\in \mathcal{N}_\e$, $\Phi_\e|_{\hat{E}_\e(z)}$ admits a unique maximum point which occurs precisely at $z$;

\item [2)] for any $z\in E_\e\setminus E_\e^-$, the set $\hat{E}_\e(z)$ intersects $\mathcal{N}_\e$ at exactly one point $\hat{m}_\e(z)$, which is the unique global maximum point of $\Phi_\e|_{\hat{E}_\e(z)}$.
\end{itemize}
\el
\subsection{Lower and upper bounds for $c_\e$}
\bo\lab{co51} There exists $c_0>0$ (independent of $\e$) such that for $\e>0$ sufficiently small,
$$c_\e=\inf_{z\in E_\e\setminus E_\e^-}\max_{\omega\in\hat{E}_\e(z)}\Phi_\e(\omega)\in(c_0,4\pi/\al_0).$$
\eo
\bp
The min-max characterization is standard and we refer to \cite{DJJ}. Here we are concerned with estimating form below and above the critical level $C_\e$. 

\noindent {\it Lower bound.} On one hand, for any $z\in E_\e$, we know $\hat{E}_\e(z)=\hat{E}_\e(z^+)$. Then, for any $a>0$
\begin{align*}
c_\e&=\inf_{z\in E_\e\setminus E_\e^-}\max_{\omega\in\hat{E}_\e(z)}\Phi_\e(\omega)=\inf_{z\in E_\e^+\setminus\{0\}}\max_{\omega\in\hat{E}_\e(z)}\Phi_\e(\omega)\\
&=\inf_{z\in S_{a,\e}^+}\max_{\omega\in\hat{E}_\e(z)}\Phi_\e(\omega)\ge\inf_{z\in S_{a,\e}^+}\max_{\omega\in\R^+ z}\Phi_\e(\omega),
\end{align*}
where
$
S_{a,\e}^+:=\{z\in E_\e^+: \|z\|_\e=a\}.
$
On the other hand, recalling that $f,g$ have critical growth with critical exponent $\al_0$, by $(H1)$, for some $\al'>\al_0$, there exists $C>0$ such that
\be\lab{fg}
F(t)\le \frac{1}{4}V_0|t|^2+C |t|^3\left(e^{\al't^2}-1\right),\, G(t)\le \frac{1}{4}V_0|t|^2+C |t|^3\left(e^{\al't^2}-1\right), t\in\R.
\ee
By the Pohozaev-Trudinger-Moser inequality, there exists $a>0$ sufficiently small such that
$$
\int_{\R^2}\left(e^{2\al'u^2}-1\right)\le1,
$$
for any $u\in H^1(\R^2)$ with $\|u\|_{H^1}\le a$. Then, for any $z=(u,u)\in S_{a,\e}^+$,
{\allowdisplaybreaks
\begin{align*}
&\max_{\omega\in\R^+ z}\Phi_\e(\omega)\ge\Phi_\e(z)=\int_{\R^2}|\na u|^2+V_\e(x)u^2-\int_{\R^2}F(u)+G(u)\\
&\ge\|u\|_{1,\e}^2-V_0/2\int_\Omega u^2-2C\int_{\R^2} |u|^3\left(e^{\al'u^2}-1\right)\\
&\ge C'\|u\|_{1,\e}^2-2C\left(\int_\Omega u^6\right)^{1/2}\ge\|u\|_{1,\e}^2(C'-2C C_6^3\|u\|_{1,\e}),
\end{align*}
}%
where $C'=\min\{1,V_0\}/2$ and $C_6$ is the Sobolev's constant of the embedding $H^1(\R^2)\hookrightarrow L^6(\R^2)$. Thus, taking $a>0$ fixed but small enough, for any $z=(u,u)\in S_{a,\e}^+$, we have $\|u\|_{1,\e}^2=a^2/2$ and
$$
\max_{\omega\in\R^+ z}\Phi_\e(\omega)\ge\|u\|_{1,\e}^2\left[C'-2CC_6^3\|u\|_{1,\e}\right]\ge a^2/6>0.
$$
Thus, for any $\e>0$, $c_\e\ge c_0=a^2/6$.

\noindent {\it Upper bound.} By $(H5)$ and $V(0)=V_0$, for some fixed $r>0$ and $\e_0>0$ such that
\begin{equation}
\label{ChoiceOfBeta}
\beta_0>\frac{4e^{\frac{r^2}{2}\max_{|x|\le\e r}V(x)}}{\al_0 r^2},\,\,\e\in(0,\e_0),
\end{equation}
we consider the following so-called Moser sequence
\begin{equation}
\omega_k(x) = \frac{1}{\sqrt{2\pi}}\left\{
\begin{array}{ll}
(\log k)^{1/2},  & \mid x \mid \leq r/k;\\
\frac{ \log \frac{r}{\mid x \mid} }{(\log k)^{1/2}}, & r/k \leq
\mid x \mid \leq r; \\
0, & \mid x \mid \geq r.
\end{array}
\right.
\end{equation}
Then, one easily checks that $\|\na\omega_k\|_2=1$ and $\|\omega_k\|_2^2=r^2/(4\log{k})+o(r^2/\log{k})$. Let $d_k(r):= r^2/4 + o_k(1)$ where $o_k(1) \to 0$, as $k \to + \infty$ and $\ti{\omega}_{k,\e}:=\omega_k/\|\omega_k\|_{1,\e}$, then $\|\ti{\omega}_{k,\e}\|_{1,\e}=1$ and for $k$ large enough,
\begin{equation}
\label{EstOfwn}
\ti{\omega}_{k,\e}^2(x)\geq \frac 1{2 \pi} \: \Bigl( \log k -  \: d_{k,\e}(r) \Bigr) \quad \text{for } |x| \leq \frac r k,
\end{equation}
where $d_{k,\e}(r)=d_k(r)\max_{|x|\le\e r}V(x)\ge V_0d_k(r)$.

\noindent Suppose by contradiction that for some fixed $\e\in(0,\e_0)$ and for all $ k $,
$$
\sup_{z\in\hat{E}((\ti{\omega}_{k,\e},\ti{\omega}_{k,\e}))}\Phi_\e(z) \geq 4\pi/\alpha_0.
$$
Then $\Phi_\e(\hat{m}((\ti{\omega}_{k,\e},\ti{\omega}_{k,\e})))\geq 4\pi/\alpha_0$ for all $ k $, where $\hat{m}((\ti{\omega}_{k,\e},\ti{\omega}_{k,\e}))\in\mathcal{N}_\e$ and
$$
\hat{m}((\ti{\omega}_{k,\e},\ti{\omega}_{k,\e})) = \tau_k(\ti{\omega}_{k,\e},\ti{\omega}_{k,\e}) + (u_k,-u_k)\in\hat{E}((\ti{\omega}_{k,\e},\ti{\omega}_{k,\e})).
$$
Namely,
\begin{equation}\label{Itatiaia}
\tau^2_{k}  - \int_{\R^2}(|\nabla u_{k}| ^2+V_\e(x)u_k^2) -
 \int_{\R^2} [F(\tau_{k}\ti{\omega}_{k,\e} + u_{k})+G(\tau_{k}\ti{\omega}_{k,\e}-u_{k})] \geq 4\pi/\alpha_0
\end{equation}
and
\begin{equation}\label{Itatuba}
\tau^2_{k} -
  \int_{\R^2}(|\nabla u_{k}| ^2+V_\e(x)u_k^2) =
  \int_{\R^2} [f(\tau_{k}\ti{\omega}_{k,\e} + u_{k})(\tau_{k}\ti{\omega}_{k,\e} + u_{k})
 + g(\tau_{k}\ti{\omega}_{k,\e} -u_{k})(\tau_{k}\ti{\omega}_{k,\e} -u_{k})].
\end{equation}

\noindent {\it Claim:} $\lim_{k\rg\iy}\tau_k=4 \pi / \alpha_0 $. Indeed, from (\ref{Itatiaia}), we get $\tau_{k}^2 \geq 4 \pi / \alpha_0 $.
From $(H5)$, given $ \rho>0 $, there exists $ R_\rho $ such
that
\[
tf(t) \geq (\beta_0 - \rho) e^{\alpha_0 t^2}
\mbox{ for all } t \geq R_\rho .
\]
and the same holds true also for $tg(t)$. 
Noting that $$\tau_k \ti{\omega}_{k,\e}=\frac{\tau_k}{\|\omega_k\|_\e}\frac{\sqrt{\log{k}}}{\sqrt{2\pi}}\rg+\iy,\ \ \mbox{as}\ \ k\rg\iy,\ \ x\in B_{r/k},$$
by choosing $k$ sufficiently large, we get 
$\max{\{\tau_{k}\ti{\omega}_{k,\e} +u_{k}, \; \tau_{k}\ti{\omega}_{k,\e} -u_{k}\}} \geq
R_\rho $ for all $ x \in B_{r/k} $. So that by \re{EstOfwn},
\begin{align}\lab{bound}
\tau_k^2&\ge\int_{B_{r/k}} [f(\tau_{k}\ti{\omega}_{k,\e} + u_{k})(\tau_{k}\ti{\omega}_{k,\e} + u_{k})
 + g(\tau_{k}\ti{\omega}_{k,\e} -u_{k})(\tau_{k}\ti{\omega}_{k,\e} -u_{k})]\nonumber\\
&\geq (\beta_0- \rho) \int_{B_{r/k}}e^{\al_0 (\tau_k\ti{\omega}_{k,\e})^2} \, \ud x\nonumber\\
&\ge\pi r^2(\beta_0 - \rho)\: e^{\frac{\al_0}{2\pi}\tau_k^2[\log k - d_{k,\e}(r)]-2\log{k}},
\end{align}
which implies that $\{\tau_k\}$ is bounded. By \re{bound}, as a consequence of the boundedness of $\{\tau_k\}$, we know $\limsup_{k\rg\iy}\tau_k^2\le4 \pi / \alpha_0$. In fact, if not we have
$$
\limsup_{k\rg\iy}e^{\frac{\al_0}{2\pi}\tau_k^2[\log k - d_{k,\e}(r)]-2\log{k}}=\iy,
$$
which is a contradiction, and the claim is proved.

\noindent As $\omega_k \to 0 $ a.e. in $\R^2$, by the Lebesgue dominated convergence theorem
$$\int_{\{x\in B_r:\tau_k\ti{\omega}_{k,\e} < R_\rho\}} \min\{f(\tau_k\ti{\omega}_{k,\e}) \tau_k\ti{\omega}_{k,\e}, g(\tau_k\ti{\omega}_{k,\e}) \tau_k\ti{\omega}_{k,\e}\}\, \ud x \rg0,\quad  k\rg\iy$$
and $$\int_{\{x\in B_r:\tau_k\ti{\omega}_{k,\e} < R_\rho\}} e^{\al_0 (\tau_k\ti{\omega}_{k,\e})^2} \, \ud x\rg\pi r^2.$$
Then, from \eqref{Itatuba} and $(H4)$ we have obtain 
\begin{align*}
\tau_k^2 & \ge\int_{B_r} [f(\tau_k\ti{\omega}_{k,\e}+u_k)(\tau_k\ti{\omega}_{k,\e}+u_k)+g(\tau_k\ti{\omega}_{k,\e}-u_k)(\tau_k\ti{\omega}_{k,\e}-u_k)] \, \ud x\\
& \geq (\beta_0- \rho) \int_{B_r} e^{\al_0 (\tau_k\ti{\omega}_{k,\e})^2} \, \ud x - (\beta_0- \rho) \int_{\{x\in B_r:\tau_k\ti{\omega}_{k,\e} < R_\rho\}} e^{\al_0 (\tau_k\ti{\omega}_{k,\e})^2} \, \ud x\\
&\ \ \ \ +\int_{\{x\in B_r:\tau_k\ti{\omega}_{k,\e} < R_\rho\}}\min\{f(\tau_k\ti{\omega}_{k,\e}) \tau_k\ti{\omega}_{k,\e}, g(\tau_k\ti{\omega}_{k,\e}) \tau_k\ti{\omega}_{k,\e}\}\, \ud x\\ &=(\beta_0- \rho)\Big[\int_{B_r} e^{\al_0 (\tau_k\ti{\omega}_{k,\e})^2} \, \ud x-\pi r^2\Big].
\end{align*}
In the following, we estimate the term $\int_{B_r} e^{\al_0 (\tau_k\ti{\omega}_{k,\e})^2} \, \ud x$. Observe first that from \eqref{EstOfwn} one has 
$$\int_{B_{r/k}}  e^{\al_0 (\tau_k\ti{\omega}_{k,\e})^2} \, \ud x \geq \pi r^2 \: e^{\frac{\al_0}{2\pi}\tau_k^2[\log{k}- d_{k,\e}(r)]-2\log{k}}.$$ Noting that $\tau_k^2\ge4\pi/\al_0$ and $\tau_k^2\rg4\pi/\al_0$, we have
$$\liminf_{k\rg\iy}\int_{B_{r/k}}  e^{\al_0 (\tau_k\ti{\omega}_{k,\e})^2}\, \ud x \geq \pi r^2e^{-\max_{|x|\le\e r}V(x)r^2/2}.$$
Secondly, by using the change of variable $s=r e^{-\|\omega_k\|_\e \sqrt{\log k} \: t}$, one has 
\begin{equation*}
\begin{split}
\int_{B_r \setminus B_{r/k}} e^{4 \pi (\ti{\omega}_{k,\e})^2} \, \ud x & = 2 \pi r^2 \|\omega_k\|_\e \sqrt{\log k} \: \int_{0}^{\frac{\sqrt{\log k}}{\|\omega_k\|_\e}} \: e^{2(\:t^2 - \|\omega_k\|_\e \sqrt{\log k} \: t \:)} \, \ud t
\\
& \geq 2 \pi r^2 \|\omega_k\|_\e \sqrt{\log k} \: \int_{0}^{\frac{\sqrt{\log k}}{\|\omega_k\|_\e}} \: e^{-2\|\omega_k\|_\e \sqrt{\log k} \: t} \, \ud t \\
& = \pi r^2 \bigl( 1 - e^{-2 \log k}\bigr).
\end{split}
\end{equation*}
Thus $$\liminf_{k\rg\iy}\int_{B_r} e^{\al_0(\tau_k \ti{\omega}_{k,\e})^2} \, \ud x \geq \pi r^2(e^{-\max_{|x|\le\e r}V(x)r^2/2}+1),$$
which implies 
$$4\pi/\al_0= \lim_{k \to + \infty} \tau_k^2 \geq (\beta_0 - \rho) \pi r^2e^{-\max_{|x|\le\e r}V(x)r^2/2}.$$
As $\rho$ is arbitrary, we have
$$\beta_0 \leq \frac{4e^{\frac{r^2}{2}\max_{|x|\le\e r}V(x)}}{\al_0 r^2},$$
which contradicts \eqref{ChoiceOfBeta}. Therefore, $c_\e<4\pi/\al_0$ for $\e\in(0,\e_0)$.
\ep

\subsection{Existence of solutions to system \re{q51}}\

\noindent Let us define
$$
\hat{m}_\e: z\in E_\e\setminus E_\e^-\mapsto\hat{m}_\e(z)\in\hat{E}_\e(z)\cap\mathcal{N}_\e.
$$
\bl\lab{l5.4}
There exists $\dd>0$ (independent of $\e$) such that $\|z^+\|_\e\ge\dd$ for all $z\in\mathcal{N}_\e$. In particular,
$$
\|\hat{m}_\e(z)^+\|_\e\ge\dd\ \ \ \mbox{for all}\ \ z\in E_\e\setminus E_\e^-.
$$
Moreover, for each compact subset $\mathcal{W}\subset E_\e\setminus E_\e^-$, there exists a constant $C_{\mathcal{W},\e}>0$ such that
$$
\|\hat{m}_\e(z)\|_\e\le C_{\mathcal{W},\e}\ \ \ \mbox{for all}\ \ z\in\mathcal{W}.
$$
\el
\noindent Let
$$
S_\e^+:=\{z\in E_\e^+: \|z\|_\e=1\},
$$
then $S_\e^+$ is a $C^1$-submanifold of $E_\e^+$ and the tangent manifold of $S_\e^+$ at $z\in S_\e^+$ is
$$
T_z(S_\e^+)=\{\omega\in E_\e^+: (\omega,z)_\e=0\}.
$$
Let
$$
m_\e:=\hat{m}_\e|_{S_\e^+}: S_\e^+\longrightarrow \mathcal{N}_\e,
$$
then by Lemma \ref{l5.4}, $\hat{m}_\e$ is continuous and $m_\e$ is a homeomorphism between $S_\e^+$ and $\mathcal{N}_\e$. Define
$$
\Psi_\e: S_\e^+\longrightarrow\R,\quad  \Psi_\e(z):=\Phi_\e(m_\e(z)),\: z\in S_\e^+,
$$
then, as a consequence of \cite[Corollary 4.3]{Weth}, for any fixed $\e>0$, we have the following 
\bo\lab{p5.5}\noindent
\begin{itemize}

\item [1)] $\Psi_\e\in C^1(S_\e^+,\R)$ and
$$
\lan\Psi_\e'(z),\omega\ran_\e=\|m_\e(z)^+\|\lan\Phi_\e'(m_\e(z)),\omega\ran_\e\ \ \mbox{for all}\ \ \omega\in T_z(S_\e^+);
$$

\item [2)] If $\{\omega_n\}\subset S_\e^+$ is a Palais-Smale sequence for $\Psi_\e$, then $\{m_\e(\omega_n)\}\subset \mathcal{N}_\e$ is a Palais-Smale sequence for $\Phi_\e$. Namely, if $\Psi_\e(\omega_n)\rg d$ for some $d>0$ and $\|\Psi_\e'(\omega_n)\|_\ast\rg 0$, as $n\rg\iy$, then $\Phi_\e(m_\e(\omega_n))\rg d$ and $\|\Phi_\e'(m_\e(\omega_n))\|\rg0$, as $n\rg\iy$, where
    $$
    \|\Psi_\e'(\omega_n)\|_\ast=\sup_{\stackrel{\phi\in T_{\omega_n}(S_\e^+)}{\|\phi\|_\e=1}}\lan\Psi_\e'(\omega_n),\phi\ran_\e\ \ \mbox{and}\ \  \|\Phi_\e'(m_\e(\omega_n))\|=\sup_{\stackrel{\phi\in E_\e}{\|\phi\|_\e=1}}\lan\Phi_\e'(m_\e(\omega_n)),\phi\ran_\e;
    $$

\item [3)] $\omega\in S_\e^+$ is a critical point of $\Psi_\e$ if and only if $m_\e(\omega)\in \mathcal{N}_\e$ is a critical point of  $\Phi_\e$;

\item [4)] $\inf_{S_\e^+}\Psi_\e=\inf_{\mathcal{N}_\e}\Phi_\e$.
\end{itemize}
\eo
\noindent Since $S_\e^+$ is a regular $C^1$-submanifold of $E_\e^+$, by Proposition \ref{co51} and Proposition \ref{p5.5}, it follows from the Ekeland variational principle (see \cite[Theorem 3.1]{E}) that there exists $\{\omega_n\}\subset S_\e^+$ such that
$$
\Psi_\e(w_n)\rg c_\e>0 \ \ \mbox{and}\ \ \|\Psi_\e'(\omega_n)\|_\ast\rg 0,\ \ \mbox{as}\ \ n\rg\iy.
$$
Let $z_n=m(\omega_n)\in\mathcal{N}_\e$, then
\be
\Phi_\e(z_n)\rg c_\e>0 \ \ \mbox{and}\ \ \|\Phi_\e'(z_n)\|\rg 0,\ \ \mbox{as}\ \ n\rg\iy.
\ee
Similar as in \cite{DJJ}, one has the following two propositions:
\bo\lab{o54} There exists $C$ (independent of $\e$) such that for all $\e>0$ and $n\in \mathbb{N}$:
\begin{itemize}

\item [1)] $\|z_n\|_\e=\|(u_n,v_n)\|_\e\le C(1+c_\e)$;

\item [2)] $\int_{\R^2}f(u_n)u_n\, \ud x\le C(1+c_\e)$ and $\int_{\R^2}g(v_n)v_n\, \ud x\le C(1+c_\e)$ ;

\item [3)] $\int_{\R^2}F(u_n)\, \ud x\le C(1+c_\e)$ and $\int_{\R^2}G(v_n)\, \ud x\le C(1+c_\e)$.
\end{itemize}
\eo
\noindent Up to a subsequence, there exists $z_\e=(u_\e,v_\e)\in E_\e$ such that $z_n\rightharpoonup z_\e$ in $E_\e$ and $z_n\xrightarrow{a.e.}z_\e$ in $\R^2$, as $n\rg\iy$, which is actually a weak solution to \re{q51}, precisely we have 
\bo\lab{o51}
The weak limit $z_\e$ is a critical point of $\Phi_\e$.
\eo

\subsection{Asymptotic behavior of $c_\e$}

By Proposition \ref{o51}, it suffices to show $z_\e\not\equiv0$. For this purpose, in the following, we investigate the relation between $c_\ast$ and $c_\e$, where $c_\ast, c_\e$ are the corresponding least energies to System \re{q11} and \re{q51} respectively.
\bl\lab{l5.9}
With the assumptions of Theorem \ref{Th1}, we have $$\limsup_{\e\rg0}c_\e\le c_\ast.$$
\el
\bp
By Theorem \ref{a}, there exists $z=(u,v)\in\mathcal{N}$ such that $$c_\ast=\max_{\omega\in\hat{E}(z)}\Phi(\omega)=\max_{\omega\in\hat{E}(z^+)}\Phi(\omega).$$ Noting that $z\in E\setminus E^-$, we know for any $\e>0$, $z\in E_\e\setminus E_\e^-$. Then, by Lemma \ref{l5.1}, for any $\e>0$
\begin{align*}
c_\e\le\max_{\omega\in\hat{E}_\e(z)}\Phi_\e(\omega)=\Phi_\e(\hat{m}_\e(z)).
\end{align*}
Recalling that $\hat{m}_\e(z)\in\hat{E}_\e(z)\cap\mathcal{N}_\e$, there exist $s_\e\ge0$, $t_\e\in\R$ and $\vp_\e\in H_\e$, $\|\vp_\e\|_\e=1$ such that $\hat{m}_\e(z)=s_\e z+t_\e(\vp_\e,-\vp_\e)$.\vskip0.1in
{\bf Step 1.} We borrow some ideas from \cite{Ramos1} to prove that $t_\e,s_\e$ are bounded for $\e>0$ sufficiently small. We proceed by contradiction and distinguish between two cases.

\noindent {\it Case I.} Both $s_\e,t_\e$ are unbounded for $\e$ small. If $|t_\e|/s_\e\rg\iy$, as $\e\rg0$, then
\begin{align*}
c_\e&\le\Phi_\e(s_\e z+t_\e(\vp_\e,-\vp_\e))\\
&=s_\e^2\|z\|_\e^2-t_\e^2+t_\e s_\e O(1)-\int_{\R^2}F(s_\e u+t_\e\vp_\e)+G(s_\e v-t_\e\vp_\e)\\
&\le s_\e^2\|z\|_\e^2-t_\e^2+t_\e s_\e O(1)=s_\e^2(O(1)-1)\rg-\iy,
\end{align*}
which contradict the fact $c_\e\ge c_0>0$. If $|t_\e|/s_\e\rg0$, as $\e\rg0$, then
\begin{align*}
c_\e&\le s_\e^2\|z\|_\e^2-t_\e^2+t_\e s_\e O(1)-\int_{\R^2}F(s_\e u+t_\e\vp_\e)+G(s_\e v-t_\e\vp_\e)\\
&=s_\e^3\left(o(1)-\int_{\R^2}\frac{F(s_\e u+t_\e\vp_\e)}{s_\e^3}+\frac{G(s_\e v-t_\e\vp_\e)}{s_\e^3}\right).
\end{align*}
Since $c_\e\ge c_0>0$, as $\e\rg0$ we have 
$$
\int_{\R^2}\frac{F(s_\e u+t_\e\vp_\e)}{s_\e^3}\rg0,\,\,\frac{G(s_\e v-t_\e\vp_\e)}{s_\e^3}\rg0.
$$
Recalling that $f$ has Moser critical growth at infinity, there exists $C>0$ such that $|F(t)|\ge C|t|^3$ for $|t|\ge1$. Let $A_\e:=\{x\in\R^2: |s_\e u(x)+t_\e\vp_\e(x)|\ge1\}$, then
$$
\int_{A_\e}\frac{F(s_\e u+t_\e\vp_\e)}{s_\e^3}\ge C\int_{A_\e}\left|u(x)+\frac{t_\e}{s_\e}\vp_\e(x)\right|^3,
$$
where the left hand side vanishes as $k\rg\iy$, which yields $\lim_{\e\rg0}\int_{A_\e}|u(x)|^3=0$. At the same time,
\begin{align*}
&\int_{\R^2\setminus A_\e}|u(x)|^3\le\int_{\R^2\setminus A_\e}u^2(x)\left(\frac{1}{s_\e}+\frac{|t_\e|}{s_\e}|\vp_\e|\right)\\
&\le\frac{1}{s_\e}\int_{\R^2}u^2(x)+\frac{|t_\e|}{s_\e}\left(\int_{\R^2}u^4(x)\right)^{1/2}\left(\int_{\R^2}\vp^2_\e(x)\right)^{1/2}\rg0,\hbox{ as $\e\rg0$.}
\end{align*}
Hence $\int_{\R^2}|u|^3=0$ and in turn $u\equiv0$. Similarly, $v\equiv0$. So that we get $c_\ast=0$, which is a contradiction. If $|t_\e|/s_\e\rg l>0$, as $\e\rg0$, then following the same line as above,
$$\int_{A_\e}\left|u(x)+\frac{t_\e}{s_\e}\vp_\e(x)\right|^3\rg0.$$
Moreover,
\begin{align*}
\int_{\R^2\setminus A_\e}\left|u(x)+\frac{t_\e}{s_\e}\vp_\e(x)\right|^3\le\frac{1}{s_\e}\int_{\R^2\setminus A_\e}\left|u(x)+\frac{t_\e}{s_\e}\vp_\e(x)\right|^2\rg0,\hbox{ as $\e\rg0$.}
\end{align*}
Then
$$
\int_{\R^2}\left|u(x)+\frac{t_\e}{s_\e}\vp_\e(x)\right|^3\rg0,\,\,\e\rg0
$$
and analogously
$$
\int_{\R^2}\left|v(x)-\frac{t_\e}{s_\e}\vp_\e(x)\right|^3\rg0,\,\,\e\rg0.
$$
So we get $\int_{\R^2}|u+v|^3=0$, that is $u=-v$. This implies $z=(u,v)\in E^-$ which contradicts the fact $z\in\mathcal{N}$.

\noindent {\it Case II.} Just one between $s_\e$ and $t_\e$ stays bounded for $\e$ small. If $|t_\e|/s_\e\rg\iy$, as $\e\rg0$, then $|t_\e|\rg\iy$, as $\e\rg0$ and as above one has 
\begin{align*}
c_\e\le s_\e^2\|z\|_\e^2-t_\e^2+t_\e s_\e O(1)=t_\e^2(O(1)-1)\rg-\iy,
\end{align*}
which contradicts the fact $c_\e\ge c_0>0$. If $|t_\e|/s_\e$ is bounded for $\e$ small, then $s_\e\rg\iy$ and $|t_\e|/s_\e\rg0$, as $\e\rg0$. Reasoning as in {\it Case I}, we get $u=v=0$ and $c_\ast=0$, which is again a contradiction.
\vskip0.1in
{\bf Step 2.} Recall that
\begin{align*}
c_\e\le\max_{\omega\in\hat{E}_\e(z)}\Phi_\e(\omega)=\Phi_\e(\hat{m}_\e(z))
\end{align*}
where $\hat{m}_\e(z)=s_\e z+t_\e(\vp_\e,-\vp_\e)$. Then
\begin{align*}
c_\e\le&\Phi_\e(s_\e z+t_\e(\vp_\e,-\vp_\e))=\Phi(s_\e z+t_\e(\vp_\e,-\vp_\e))\\
&+\int_{\R^2}(V_\e(x)-1)(s_\e u+t_\e\vp_\e)(s_\e v-t_\e\vp_\e)\\
\le&\max_{\omega\in\hat{E}(z)}\Phi(\omega)+I_\e=c_\ast+I_\e,
\end{align*}
where $I_\e:=\int_{\R^2}(V_\e(x)-1)(s_\e u+t_\e\vp_\e)(s_\e v-t_\e\vp_\e)$. Since $0\in\mathcal{M}$, by Lebesgue's dominated convergence theorem and Step 1, we have
\begin{align*}
I_\e&=\int_{\R^2}(V_\e(x)-1)[s_\e^2uv-t_\e^2\vp_\e^2+t_\e s_\e (v-u)\vp_\e]\\
&\le\int_{\R^2}(V_\e(x)-1)[s_\e^2uv+t_\e s_\e (v-u)\vp_\e]\\
&\le s_\e^2\int_{\R^2}(V_\e(x)-1)uv+|t_\e s_\e|\left(\int_{\R^2}|V_\e(x)-1|^2(v-u)^2\right)^{1/2}\rg0\,\,\hbox{ as $\e\rg0$.}
\end{align*}
Therefore, $\limsup_{\e\rg0}c_\e\le c_\ast$.
\ep
\subsection{Existence of ground state solutions for \re{q51}}

For any $\la>0$, let us consider the following problem in $\mathbb{R}^2$
\be\lab{q110} \left\{
\begin{array}{ll}
-\DD u+\la u=g(v)\\
-\DD v+\la v=f(u)
\end{array}
\right. \ee
whose corresponding energy functional is
$$
\Phi_\la(z):=\int_{\R^2}\na u\na v+\la uv-I(z),\ \ z=(u,v)\in E.
$$
As above one can define the generalized Nehari Manifold $\mathcal{N}_\la$ and the least energy
$$
c_\la:=\inf_{z\in\mathcal{N}_\la}\Phi_\la(z).
$$
Moreover, with the same assumptions of Theorem \ref{Th2}, if $c_\la\in(0,4\pi/\al_0)$ for some $\la>0$, then there exists $z_\la=(u_\la,v_\la)\in\mathcal{N}_\la$ such that $\Phi_\la(z_\la)=c_\la$.
\bl\lab{bj1}
With the assumptions of Theorem \ref{Th2}, for any $\la>0$ the map $\lambda\mapsto c_\la\in(0,4\pi/\al_0)$ is strictly increasing.
\el
\bp
For any $\la>0$ with $c_\la\in(0,4\pi/\al_0)$, let $z_\la=(u_\la,v_\la)$ be a solution of \re{q110}, then $\ti{z}_\la=(\ti{u}_\la,\ti{v}_\la)=(u_\la(\cdot/\sqrt{\la},v_\la(\cdot/\sqrt{\la}))$ satisfies in the whole plane the following system 
\be\lab{q111} \left\{
\begin{array}{ll}
-\DD \ti{u}_\la+\ti{u}_\la=\la^{-1}g(\ti{v}_\la)\\
-\DD \ti{v}_\la+\ti{v}_\la=\la^{-1}f(\ti{u}_\la)
\end{array}
\right. \ee
whose corresponding energy functional is
$$
\ti{\Phi}_\la(\ti{z}_\la):=\int_{\R^2}\na\ti{u}_\la\na\ti{v}_\la+\ti{u}_\la\ti{v}_\la-\la^{-1}I(\ti{z}_\la).
$$
Similar as above, we can define the generalized Nehari Manifold $\mathcal{\ti{N}}_\la$ and the least energy
$$
\ti{c}_\la:=\inf_{z\in\mathcal{\ti{N}}_\la}\ti{\Phi}_\la(z).
$$
We have $c_\la=\ti{c}_\la\in(0,4\pi/\al_0)$. Then \re{q111} admits a ground state solution $\ti{z}_\la=(\ti{u}_\la,\ti{v}_\la)$. Moreover,
$$
\ti{c}_\la:=\inf_{z\in E\setminus E^-}\max_{\omega\in\hat{E}(z)}\ti{\Phi}_\la(\omega)=\max_{\omega\in\hat{E}(\ti{z}_\la)}\ti{\Phi}_\la(\omega).
$$
To show that $c_\la$ is strictly increasing, it is enough to prove that $\ti{c}_\la$ is strictly increasing. For any $0<\mu<\la$, the set $\hat{E}(\ti{z}_\la)$ intersects $\mathcal{\ti{N}_\mu}$ at exactly one point $\hat{m}_\mu(z)$, which is the unique global maximum point of $\ti{\Phi}_\mu|_{\hat{E}(\ti{z}_\la)}$. Since $F(s), G(s)>0$ for any $s\not=0$,
\begin{align*}
\ti{c}_\mu&\le\max_{\omega\in\hat{E}(\ti{z}_\la)}\ti{\Phi}_\mu(\omega)=\ti{\Phi}_\mu(\hat{m}_\mu(z))\\
&<\ti{\Phi}_\la(\hat{m}_\mu(z))\le\max_{\omega\in\hat{E}(\ti{z}_\la)}\ti{\Phi}_\la(\omega)=\ti{c}_\la.
\end{align*}
Therefore, $c_\mu<c_\la$.
\ep
\noindent Now, we are set to prove that the weak limit obtained in Proposition \ref{o51} is non trivial, precisely
\bl\lab{bj2}
$z_\e\not\equiv0$ provided $\e>0$ is sufficiently small.
\el
\bp
Assume by contradiction that $z_\e=0$ for $\e>0$ small, then $z_n=(u_n,v_n)\rightharpoonup 0$ in $E_\e$ and $z_n\xrightarrow{a.e.}0$ in $\R^2$, as $n\rg\iy$. It is well known that $\{z_n\}$ satisfies just one of the following alternatives:
\begin{itemize}

\item [1)] (Vanishing)
$$
\lim_{n\rg\iy}\sup_{y\in\R^2}\int_{B_R(y)}(u_n^2+v_n^2)\, \ud x=0\ \ \mbox{for all}\ \ R>0;
$$

\item [2)] (Nonvanishing) there exist $\nu>0$, $R_0>0$ and $\{y_n\}\subset\R^2$ such that
$$
\lim_{n\rg\iy}\int_{B_{R_0}(y_n)}(u_n^2+v_n^2)\, \ud x\ge\nu.
$$
\end{itemize}
Due to $c_\e\in(c_0,4\pi/\al_0)$ we can rule out {\it Vanishing}. So that {\it Nonvanishing} occurs. Let $\ti{u}_n(\cdot):=u_n(\cdot+y_n)$ and $\ti{v}_n(\cdot):=v_n(\cdot+y_n)$, then $|y_n|\rg\iy$, as $n\rg\iy$ and
\be\lab{y522}
\lim_{n\rg\iy}\int_{B_{R_0}(0)}(\ti{u}_n^2+\ti{v}_n^2)\, \ud x\ge\nu.
\ee
Let $\ti{z}_n=(\ti{u}_n,\ti{v}_n)$, $\{\ti{z}_n\}$ is bounded in $E$. Up to a subsequence, by \re{y522} we assume that $\ti{z}_n\rg \ti{z}\not=0$ weakly in $E$ for some $\ti{z}=(\ti{u},\ti{v})\in E$ and $\Phi_{V_\iy}'(\ti{z})=0$, where
$$
\Phi_{V_\iy}(z)=\int_{\R^2}\na u\na v+V_\iy uv-I(z),\,\, z=(u,v)\in E.
$$
By $(H2)$ and Fatou's Lemma, for fixed $\e>0$,
\begin{align*}
c_\e+o_n(1)&=\Phi_\e(\ti{z}_n)-\frac{1}{2}\lan \Phi_\e'(\ti{z}_n),\ti{z}_n\ran\\
&=\int_{\R^2}\frac{1}{2}f(\ti{u}_n)\ti{u}_n-F(\ti{u}_n)+\int_{\R^2}\frac{1}{2}g(\ti{v}_n)\ti{v}_n-G(\ti{v}_n)\\
&\ge\int_{\R^2}\frac{1}{2}f(\ti{u})\ti{u}-F(\ti{u})+\int_{\R^2}\frac{1}{2}g(\ti{v})\ti{v}-G(\ti{v})+o_n(1)\\
&=\Phi_{V_\iy}(\ti{z})-\frac{1}{2}\lan \Phi_{V_\iy}'(\ti{z}),\ti{z}\ran+o_n(1)\ge c_{V_\iy}+o_n(1).
\end{align*}
It follows that $c_\e\ge c_{V_\iy}$ for $\e>0$ small enough. By Lemma \ref{l5.9} and Lemma \ref{bj1}, we get $c_{V_\iy}>c_\ast$. Again by Lemma \ref{l5.9} we get a contradiction.
\ep
\noindent By virtue of Lemma \ref{bj2} we get straightforward the following 
\bc
For $\e>0$ small enough, $\Phi_\e(z_\e)=c_\e$, namely $z_\e$ is a ground state solution of \re{q51}.
\ec
\subsection{Concentration}
Reasoning as in Proposition \ref{bo1} we have
\bo\lab{boc1}
Let $\e>0$ and $z_\e=(u_\e,v_\e)$ be a ground state solution to \re{q51}. Then, $u_\e, v_\e\in L^{\iy}(\R^2)\cap C_{loc}^{1,\g}(\R^2)$ for some $\g\in(0,1)$. {Moreover, $u_\e(x), v_\e(x)\rg 0$, as $|x|\rg\iy$.}
\eo
\noindent By Proposition \ref{boc1}, there exists $y_\e\in\R^2$ such that $$|u_\e(y_\e)|+|v_\e(y_\e)|=\max_{x\in\R^2}(|u_\e(x)|+|v_\e(x)|).$$ Moreover, $x_\e:=\e y_\e$ is a maximum point of $|\vp_\e(x)|+|\psi_\e(x)|$, where $(\vp_\e(\cdot),\psi_\e(\cdot))=(u_\e(\cdot/\e),v_\e(\cdot/\e))$ is a ground state solution of the original problem \re{q1}. We conclude the proof of Theorem \ref{Th1} by proving Proposition \ref{boc2}, \ref{boc3} and \ref{boc4} below.
{\bo\lab{boc2}\
\begin{itemize}
\item [1)] $\lim_{\e\rg 0}\mbox{dist}(x_\e,\mathcal{M})=0$;
\item [2)] $(u_\e(\cdot+x_\e/\e), v_\e(\cdot+x_\e/\e))$ converges (up to a subsequence) to a ground state solution of
\begin{align}\lab{luv}
\left\{
\begin{array}{ll}
-\DD u+V_0u=g(v)&\\
&\text{ in } \mathbb{R}^2\\
-\DD v+V_0v=f(u)
\end{array}
\right.
\end{align}
\item [3)] $u_\e(x+x_\e/\e), v_\e(x+x_\e/\e)\rg 0$, uniformly as $|x|\rg\iy$, for $\e>0$ sufficiently small.
\end{itemize}
\eo}
\bp
\noindent By virtue of Proposition \ref{o54} and Fatou's Lemma, there exists $C>0$ (independent of $\e$) such that $\|(u_\e,v_\e)\|_\e\le C$ for all $\e\in(0,\e_0)$.
Up to a subsequence, we may assume $z_\e=(u_\e, v_\e)\rightharpoonup z_0=(u_0, v_0)$ in $E$ and $(u_\e, v_\e)\xrightarrow{a.e.}(u_0, v_0)$ in $\R^2$, as $\e\rg0$. Due to $c_\e\in(c_0,4\pi/\al_0)$ for $\e>0$ sufficiently small, as in Lemma \ref{bj2}, we have $u_0\not\equiv0, v_0\not\equiv0$. Moreover, $\Phi'(z_0)=0$. By $(H2)$ and Fatou's Lemma,
\begin{align*}
c_\e&=\Phi_\e(z_\e)-\frac{1}{2}\lan \Phi_\e'(z_\e),z_\e\ran\\
&=\int_{\R^2}\frac{1}{2}f(u_\e)u_\e-F(u_\e)+\int_{\R^2}\frac{1}{2}g(v_\e)v_\e-G(v_\e)\\
&\ge\int_{\R^2}\frac{1}{2}f(u_0)u_0-F(u_0)+\int_{\R^2}\frac{1}{2}g(v_0)v_0-G(v_0)+o_\e(1)\\
&=\Phi(z_0)-\frac{1}{2}\lan \Phi'(z_0),z_0\ran+o_\e(1)\ge c_\ast+o_\e(1).
\end{align*}
Thanks to Lemma \ref{l5.9}, $\Phi(z_0)=c_\ast$, namely $(u_0, v_0)$ is a ground state solution of \re{luv}. Thanks to Fatou's Lemma again,
$$
\lim_{\e\rg0}\int_{\R^2}\frac{1}{2}f(u_\e)u_\e-F(u_\e)=\int_{\R^2}\frac{1}{2}f(u_0)u_0-F(u_0)
$$
and
$$
\lim_{\e\rg0}\int_{\R^2}\frac{1}{2}g(v_\e)v_\e-G(v_\e)=\int_{\R^2}\frac{1}{2}g(v_0)v_0-G(v_0).
$$
Repeating the argument in Proposition \ref{con}, we get $\|u_\e\|_\e\rg\|u_0\|_{H^1}$ and $\|v_\e\|_\e\rg\|v_0\|_{H^1}$, as $\e\rg0$. This implies $(u_\e,v_\e)\rg (u_0,v_0)$ strongly in $E$ as $\e\rg0$. Then, as in Proposition \ref{bo2} and \ref{pro_apriori}, $\{\|u_\e\|_\iy,\|v_\e\|_\iy\}$ is uniformly bounded for $\e>0$ small and $$\liminf_{\e\rg0}\min\{\|u_\e\|_\iy,\|v_\e\|_\iy\}>0.$$
As in Proposition \ref{boo4}, there exists $R_2>0$ such that
$$
\lim_{\e\rg0}\int_{B_{R_2}(x_\e/\e)}(u_\e^2+v_\e^2)\, \ud x>0.
$$
Now, we claim that $\{x_\e\}$ is bounded for $\e>0$ small enough. Suppose this does not occur, so that $|x_\e|\rg\iy$, as $\e\rg0$. Let $\bar{u}_\e(\cdot)=u_\e(\cdot+x_\e/\e)$ and $\bar{v}_\e(\cdot)=v_\e(\cdot+x_\e/\e)$ which, up to a subsequence, $(\bar{u}_\e,\bar{v}_\e)\rg \bar{z}=(\bar{u},\bar{v})$ weakly in $E$, as $\e\rg0$ and $\bar{u},\bar{v}\not\equiv0$. Moreover, $\Phi_{V_\iy}'(\bar{z})=0$. As in Lemma \ref{bj2} we get a contradiction. Therefore  $\{x_\e\}$ is bounded for $\e>0$ small. Up to a subsequence, assume $x_\e\rg x_0$, as $\e\rg0$ and let $\hat{u}_\e(\cdot)=u_\e(\cdot+x_\e/\e)$, $\hat{v}_\e(\cdot)=v_\e(\cdot+x_\e/\e)$. Then, up to a subsequence, $\hat{z}_\e=(\hat{u}_\e,\hat{v}_\e)\rg\hat{z}=(\hat{u},\hat{v})\not=0$ weakly in $E$, as $\e\rg0$ and $\Phi_{V(x_0)}'(\hat{z})=0$, where
$$
\Phi_{V(x_0)}(z)=\int_{\R^2}\na u\na v+V(x_0) uv-I(z),\,\, z=(u,v)\in E.
$$
By $(H2)$ and Fatou's Lemma,
\begin{align*}
c_\e&=\Phi_\e(z_\e)-\frac{1}{2}\lan \Phi_\e'(z_\e),z_\e\ran\\
&=\int_{\R^2}\frac{1}{2}f(\hat{u}_\e)\hat{u}_\e-F(\hat{u}_\e)+\int_{\R^2}\frac{1}{2}g(\hat{v}_\e)\hat{v}_\e-G(\hat{v}_\e)\\
&\ge\int_{\R^2}\frac{1}{2}f(\hat{u})\hat{u}-F(\hat{u})+\int_{\R^2}\frac{1}{2}g(\hat{v})\hat{v}-G(\hat{v})+o_\e(1)\\
&=\Phi_{V(x_0)}(\hat{z})-\frac{1}{2}\lan \Phi_{V(x_0)}'(\hat{z}),\hat{z}\ran+o_\e(1)\ge c_{V(x_0)}+o_\e(1).
\end{align*}
Recalling that $\limsup_{\e\rg0}c_\e\le c_\ast$, we get $c_{V(x_0)}=c_\ast$ and hence $(\hat{u},\hat{v})$ is a ground state solution of \re{luv}. Thanks to Lemma \ref{bj1}, $V(x_0)=V_0$, namely $x_0\in\mathcal{M}$ and $\lim_{\e\rg 0}\mbox{dist}(x_\e,\mathcal{M})=0$. {Moreover, $(\hat{u}_\e,\hat{v}_\e)\rg (\hat{u},\hat{v})$ strongly in $E$, as $\e\rg0$. As in Proposition \ref{bo2}, $u_\e(x+x_\e/\e), v_\e(x+x_\e/\e)\rg 0$ vanish at infinity uniformly in $\e$.}
\ep
{\bo\lab{boc3} Let $(\vp_\e,\psi_\e)$ be a ground state solution to \re{q1} and $x_\e^1, x_\e^2$ be any maximum point of $|\vp_\e|$ and $|\psi_\e|$ respectively. Then,
$$
\hbox{$\lim_{\e\rg 0}\mbox{dist}(x_\e^i,\mathcal{M})=0,\quad \lim_{\e\rg0}|x_\e^i-x_\e|=0,\quad i=1,2$.}
$$
If in addition $f$ and $g$ are odd and $(H6)$ holds, then for $\e>0$ small enough, $\vp_\e\psi_\e>0$ in $\R^2$ and
$$\lim_{\e\rg0}|x_\e^1-x_\e^2|/\e=0.$$
Moreover, for some $c,C>0$, $$|\vp_\e(x)|\le C\exp(-\frac{c}{\e}|x-x_\e^1|),\quad |\psi_\e(x)|\le C\exp(-\frac{c}{\e}|x-x_\e^2|), \,\, x\in\R^2.$$
\eo}
\bp
Note that $x_\e^1/\e,x_\e^2/\e$ are the maxima points of $u_\e,v_\e$ respectively. Thanks to the decayof $u_\e,v_\e$ and the following fact
$$
\liminf_{\e\rg0}\min\{\|u_\e\|_\iy,\|v_\e\|_\iy\}>0,
$$
we get $|x_\e^i/\e-x_\e/\e|$ is bounded for $i=1,2$ and $\e>0$ small enough. Then, $\lim_{\e\rg 0}\mbox{dist}(x_\e^i,\mathcal{M})=0\,,i=1,2$, $\lim_{\e\rg0}|x_\e^i-x_\e|=0\,,i=1,2$ and $\lim_{\e\rg0}|x_\e^1-x_\e^2|=0$ .

\noindent Next we assume {$f$ and $g$ are odd, that $(H6)$ holds}, and also that, up to a subsequence, $(x_\e^1-x_\e^2)/\e\rg y_0\in\R^2$, as $\e\rg0$.
Let $\ti{u}_\e(\cdot)=u_\e(\cdot+x_\e^1/\e)$ and $\ti{v}_\e(\cdot)=v_\e(\cdot+x_\e^2/\e)$, then $(\ti{u}_\e(\cdot),\ti{v}_\e(\cdot+(x_\e^1-x_\e^2)/\e))\rg (u,v)\not=0$ strongly in $E$ and
in $C_{loc}^1(\R^2)$, as $\e\rg0$. Moreover, $(u,v)$ is a ground state solution of \re{q11}. Without loss generality, we assume $u>0$, $v>0$
in $\R^2$. Since $0$ is a maximum point of $\ti{u}_\e$, $0$ is a maximum point also for $u$. By virtue of Theorem \ref{sign}, $0$ is the unique maximum point of
$u$ and $v$. On the other hand, up to a subsequence, $(\ti{u}_\e(\cdot+(x_\e^2-x_\e^1)/\e),\ti{v}_\e(\cdot))\rg (\ti{u},\ti{v})\not=0$ strongly in $E$ and
in $C_{loc}^1(\R^2)$, as $\e\rg0$. Then $(\ti{u}(\cdot),\ti{v}(\cdot))=(u(\cdot-y_0),v(\cdot-y_0))$, which is a ground state solution of \re{q11}.
Since $0$ is a maximum point of $\ti{v}_\e$, then $0$ is the unique maximum point of $\ti{v}$. Therefore, $y_0=0$.

\noindent Finally, we prove that $u_\e,v_\e$ do not change the sign for $\e>0$ sufficiently small. Let
$$
\bar{u}_\e=u_\e(\cdot+x_\e^1/\e),\,\,\, \bar{v}_\e=v_\e(\cdot+x_\e^1/\e),
$$
it is enough to prove $\bar{u}_\e\bar{v}_\e>0$ in $\R^2$. We assume $(\bar{u}_\e,\bar{v}_\e)\rg(u,v)\in\mathcal{S}$ strongly in $E$ and uniformly in $C_{loc}^2(\R^2)$, as $\e\rg0$ and $0$ is the unique maximum point of $u,v$. By Theorem \ref{sign}, $uv>0$ in $\R^2$. Without loss of generality, we assume $u>0$ and $v>0$ in $\R^2$. Then there exist $R>0$ and $\e_0>0$
such that $\bar{u}_\e,\bar{v}_\e>0$ in $B_R(0)$ for $\e<\e_0$. Define
$$
R_\e(\bar{u}_\e):=\sup\{r\,|\,\bar{u}_\e(x)>0,\,\, \forall\, x\in B_r(0)\},\,\,R_\e(\bar{v}_\e):=\sup\{r\,|\, \bar{v}_\e(x)>0,\,\,\forall\ x\in B_r(0)\}
$$
and $R_\e:=\min\{R_\e(\bar{u}_\e),R_\e(\bar{v}_\e)\}$, then $R_\e\ge R$ for any $\e<\e_0$. If $R_\e=\iy$ for any $\e<\e_0$, the proof is complete. Otherwise, there exists $\e_n>0$ such that $\e_n\rg0$, as $n\rg\iy$ and $R_n:=R_{\e_n}<\iy$ for any fixed $n$. Then, by the maximum principle, $R_{\e_n}(\bar{u}_{\e_n}), R_{\e_n}(\bar{v}_{\e_n})<\iy$ for any fixed $n\in\mathbb{N}$. Hence $\inf_{x\in\R^2}\bar{u}_{\e_n}(x)<0$ and $\inf_{x\in\R^2}\bar{v}_{\e_n}(x)<0$ for any $n\in\mathbb{N}$. Noting that $\bar{u}_{\e_n}(x),\bar{v}_{\e_n}(x)\rg0$, as $|x|\rg\iy$, there exist $y_n,z_n\in\R^2$ such that $\bar{u}_{\e_n}(y_n)=\min_{x\in\R^2}\bar{u}_{\e_n}(x)<0$ and $\bar{v}_{\e_n}(z_n)=\min_{x\in\R^2}\bar{v}_{\e_n}(x)<0$. Then we have
$$
g(\bar{v}_{\e_n}(y_n))\le V_0\bar{u}_{\e_n}(y_n),\,\,\,f(\bar{u}_{\e_n}(z_n))\le V_0\bar{v}_{\e_n}(z_n).
$$
By Remark \ref{remark1} we have 
$$
V_0\bar{u}_{\e_n}(y_n)\ge g(\bar{v}_{\e_n}(y_n))\ge g(\bar{v}_{\e_n}(z_n))
\ge g\left(\frac{f(\bar{u}_{\e_n}(z_n))}{V_0}\right)\ge g\left(\frac{f(\bar{u}_{\e_n}(y_n))}{V_0}\right),
$$
which yields $\inf_n|\bar{u}_{\e_n}(y_n)|>0$ by $(H1)$ . {Observe that $\bar{u}_{\e_n}(x)\rg0$, as $|x|\rg\iy$, uniformly in $\e$, and thus $\sup_n|y_n|<\iy$, namely $|y_n|<R_n$ for $n$ sufficiently large.} Hence $\bar{u}_{\e_n}(y_n)>0$, which is a contradiction. {Finally, since $u_\e,v_\e$ do not change the sign, by the standard comparison principle, we get the uniformly exponential decay at infinity.} 
\ep
In order to complete the proof of Theorem \ref{Th1} we need to prove the uniqueness of the maximum points of $\vp_\e,\psi_\e$.

\bo\lab{boc4} Let $x_\e^1, y_\e^1$ be any maxima points of $\vp_\e$. {Assume $f$ and $g$ are odd and $(H6)$ holds. Then $x_\e^1=y_\e^1$, for $\e>0$ sufficiently small.} Namely, the maximum point of $\vp_\e$ is unique. The same holds for $\psi_\e$.
\eo
\bp Let
$$
\bar{u}_\e=u_\e(\cdot+x_\e^1/\e),\,\,\, \bar{v}_\e=v_\e(\cdot+x_\e^1/\e).
$$
Then $(\bar{u}_\e,\bar{v}_\e)\rg (u,v)\in\mathcal{S}$ strongly in $E$ and uniformly in $C_{loc}^2(\R^2)$, as $\e\rg0$. Moreover, there exist $c,C>0$ such that
$$
|\bar{u}_\e(x)|\le C\exp{(-c|x-x_\e^1/\e|)},\,\,\ x\in\R^2.
$$
Hence $\|\bar{u}_\e\|_\iy\le C\exp{(-c|y_\e^1-x_\e^1|/\e)}$. As a consequence we have $$\limsup_{\e\rg0}|y_\e^1-x_\e^1|/\e<\iy.$$ 
Indeed, otherwise $\|\bar{u}_\e\|_\iy\rg0$, as $\e\rg0$, which yields 
$$
\int_{\R^2}[|\na\bar{v}_\e|^2+V_\e(x+x_\e^1/\e)|\bar{v}_\e|^2]\,\ud x=\int_{\R^2}f(\bar{u}_\e)\bar{v}_\e\,\ud x\rg0.
$$
Namely $\|v_\e\|_{1,\e}\rg0$, as $\e\rg0$ from which $\Phi_\e(u_\e,v_\e)\rg0$, as $\e\rg0$, thus a contradiction by Proposition \ref{co51}. Therefore $|y_\e^1-x_\e^1|/\e$ stays bounded for $\e>0$ small. As in Proposition \ref{boc3}, $|y_\e^1-x_\e^1|/\e\rg0$, as $\e\rg0$. Obverse that $\na\bar{u}_\e(0)=\na\bar{u}_\e((y_\e^1-x_\e^1)/\e)=0$. By Theorem \ref{sign}, $\DD u(0)<0$. Recalling that $u(x)=u(|x|)$, $u'(0)=0$ and $u''(r)<0$ for $r=|x|$ small. On the other hand, since $g\in C^1$, $\bar{u}_\e\in C^2$ and $\bar{u}_\e\rg u$ in $C_{loc}^2(\R^2)$, as $\e\rg0$, it follows from \cite[Lemma 4.2]{Ni-Takagi1} that $y_\e^1=x_\e^1$ for $\e>0$ sufficiently small.
\ep


\end{document}